\documentstyle[12pt, amstex]{amsart}

\makeatletter

\begin{document}
\bibliographystyle{amsalpha}

\newcommand{\e}{\epsilon}
\newcommand{\bc}{{\bold c}}
\newcommand{\0}{{\bold 0}}
\newcommand{\w}{{\bold w}}
\newcommand{\y}{{\bold y}}
\newcommand{\balpha}{{\boldsymbol \alpha}}
\newcommand{\bt}{{\bold t}}
\newcommand{\z}{{\bold z}}
\newcommand{\x}{{\bold x}}
\newcommand{\N}{{\bold N}}
\newcommand{\Z}{{\bold Z}}
\newcommand{\F}{{\bold F}}
\newcommand{\R}{{\bold R}}
\newcommand{\Q}{{\bold Q}}
\newcommand{\C}{{\bold C}}
\newcommand{\BP}{{\bold P}}
\newcommand{\cO}{{\mathcal O}}
\newcommand{\cX}{{\mathcal X}}
\newcommand{\cH}{{\mathcal H}}
\newcommand{\cM}{{\mathcal M}}
\newcommand{\cD}{{\mathcal D}}
\newcommand{\cB}{{\mathcal B}}
\newcommand{\sB}{{\sf B}}
\newcommand{\cY}{{\mathcal Y}}
\newcommand{\cT}{{\mathcal T}}
\newcommand{\cI}{{\mathcal I}}
\newcommand{\cS}{{\mathcal S}}
\newcommand{\sE}{{\sf E}}

\newcommand{\sA}{{\sf A}}
\newcommand{\ga}{{\sf a}}
\newcommand{\es}{{\sf s}}
\newcommand{\m}{{\bold m}}
\newcommand{\bS}{{\bold S}}
\newcommand{\ovf}{{\overline{f}}}

\newcommand{\ihra}{\stackrel{i}{\hookrightarrow}}
\newcommand\rank{\mathop{\rm rank}\nolimits}
\newcommand\im{\mathop{\rm Im}\nolimits}
\newcommand\coker{\mathop{\rm coker}\nolimits}
\newcommand\Li{\mathop{\rm Li}\nolimits}
\newcommand\Aut{\mathop{\rm Aut}\nolimits}
\newcommand\NS{\mathop{\rm NS}\nolimits}
\newcommand\Hom{\mathop{\rm Hom}\nolimits}
\newcommand\Ext{\mathop{\rm Ext}\nolimits}
\newcommand\Pic{\mathop{\rm Pic}\nolimits}
\newcommand\Spec{\mathop{\rm Spec}\nolimits}
\newcommand\Hilb{\mathop{\rm Hilb}\nolimits}
\newcommand\Ker{\mathop{\rm Ker}\nolimits}
\newcommand{\length}{\mathop{\rm length}\nolimits}
\newcommand{\res}{\mathop{\sf res}\nolimits}

\newcommand\lra{\longrightarrow}
\newcommand\ra{\rightarrow}
\newcommand\cJ{{\mathcal J}}
\newcommand\JG{J_{\Gamma}}
\newcommand{\wvskp}{\vspace{1cm}}
\newcommand{\vskp}{\vspace{5mm}}
\newcommand{\nvskp}{\vspace{1mm}}
\newcommand{\nid}{\noindent}
\newcommand{\new}{\nvskp \nid}
\newtheorem{Assumption}{Assumption}[section]
\newtheorem{Theorem}{Theorem}[section]
\newtheorem{Lemma}{Lemma}[section]
\newtheorem{Remark}{Remark}[section]
\newtheorem{Corollary}{Corollary}[section]
\newtheorem{Conjecture}{Conjecture}[section]
\newtheorem{Proposition}{Proposition}[section]
\newtheorem{Example}{Example}[section]
\newtheorem{Definition}{Definition}[section]
\newtheorem{Question}{Question}[section]
\renewcommand{\thesubsection}{\it}

\baselineskip=14pt
\addtocounter{section}{-1}

\title{Resolution of singularities 
for a family of third-order differential systems 
with small meromorphic solution spaces }

\author{Yusuke Sasano}

\address{Department of Mathematics, Faculty of Science, Kobe University, Kobe, Rokko, 657-8501, Japan}
\email{sasano@math.kobe-u.ac.jp}
\maketitle

\section{abstract}
\begin{quote}
Regarding the resolution of singularities for the differential equations of Painlev\'e type, there are important differences between the second-order Painlev\'e equations and those of higher order. Unlike the second-order case, in higher order cases there may exist some meromorphic solution spaces with codimension 2. In this paper, we will give an explicit global resolution of singularities for a 3-parameter family of third-order differential systems with meromorphic solution spaces of codimension 2.
\end{quote}

\section{Introduction}

In 1979, K. Okamoto constructed the spaces of initial conditions of Painlev\'e equations, which can be considered as the parametrized spaces of all solutions, including the meromorphic solutions [1]. They are constructed by means of successive blowing-up procedures at singular points. For Painlev\'e equations, the dimension of the space of meromorphic solutions through any singular point is always codimension one. But, in the case of higher order Painlev\'e equations, the space of meromorphic solutions through a singular point may be of codimension greater than or equal to 2 [5]. In this paper, we will give an explicit resolution of singularities for a 3-parameter family of third-order differential systems with meromorphic solution spaces of codimension 2. For second-order Painlev\'e equations, we can obtain the entire space of initial conditions by adding subvarieties of codimension 1 (equivalently, of dimension 1) to the space of initial conditions of holomorphic solutions [3]. But in the case of third-order differential equations, we need to add codimension 2 subvarieties to the space in addition to codimension 1 subvarieties. In order to resolve singularities, we need to both blow up and blow down. Moreover, to obtain a smooth variety by blowing-down, we need to resolve for a pair of singularities.

\section{ Statement of main results}

Let us consider a 3-parameter family of third-order differential system:

\begin{equation}
  \left\{
  \begin{aligned}
   \frac{dx}{dt} &=x(t-x-2z)+\alpha_1\\
   \frac{dy}{dt} &=y(-t+y+2z)+\alpha_2\\
   \frac{dz}{dt} &=z(t-2y-z)+\alpha_3.
   \end{aligned}
  \right. 
\end{equation}
Here, $x,y$ and $z$ denote unknown complex variables and $\alpha_1,\alpha_2$ and $\alpha_3$ are constant parameters. The last two equations of this system $(1)$ are equivalent to the fourth Painlev\'e equations with unknown variables $y$ and $z$. For each solution $z(t)$ of the Painlev\'e IV equations, the first equations of this system $(1)$ give Riccati equations with unknown variable $x$.

From the viewpoint of symmetry, it is interesting to point out that our differential system $(1)$ is invariant under the following transformations:

        $$ s_1: (x,y,z,t;\alpha_1,\alpha_2,\alpha_3) \rightarrow (x,y-\frac{\alpha_3}{z},z,t;\alpha_1,\alpha_2+\alpha_3,-\alpha_3),$$

        $$ s_2: (x,y,z,t;\alpha_1,\alpha_2,\alpha_3) \rightarrow (x-\frac{\alpha_2}{y},y,z+\frac{\alpha_2}{y},t;\alpha_1+\alpha_2,-\alpha_2,\alpha_3+\alpha_2).$$

In order to characterize this differential system $(1)$, we can prove the following two propositions by direct calculation.

{\bf Notation.}

$\bullet$ $W=<s_1,s_2>$,

$\bullet$ ${\cal V}:=\{(dx/dt=f_1(x,y,z,t;a_i), \ dy/dt=f_2(x,y,z,t;b_i), \ dz/dt=f_3(x,y,z,t,c_i))\},$
\noindent
where $f_i(X,Y,Z,t;e_i)=e_1X^2+e_2Y^2+e_3Z^2+e_4XY+e_5XZ+e_6YZ+e_7tX+e_8tY+e_9tZ+e_{10}$ with $e_i \in {\Bbb C} \ (1 \leq i \leq 10)$,

$\bullet$ ${\cal H} \cong {\Bbb P}^2( \subset {\Bbb P}^3)$: boundary divisor of ${\Bbb P}^3$.

\begin{Proposition}
The differential systems in $\cal V$ with symmetry under $W$ are written as follows$:$

\begin{equation}
  \left\{
  \begin{aligned}
   \frac{dx}{dt} &=a_1x^2+(a_1-c_3)z^2+a_5xz+a_6yz+a_7tx+(a_7+b_8)tz\\
                 & \hspace{9cm} +\alpha_2(b_2-a_6)+\alpha_3a_6+a_{10}\\
   \frac{dy}{dt} &=b_2y^2+(a_5-2a_1)xy+(2a_1-a_5-2c_3)yz+b_8ty+\alpha_2(a_5-2a_1)+\alpha_2c_3\\
   \frac{dz}{dt} &=c_3z^2+(2a_1-a_5)xz-2b_2yz-b_8tz+(\alpha_2-\alpha_3)b_2.
   \end{aligned}
  \right. 
\end{equation}
\end{Proposition}

\begin{Proposition}
The last two equations of this system $(2)$ are equivalent to the fourth Painlev\'e systems with unknown variables $y$ and $z$ if and only if

$$a_5=2a_1.$$
\end{Proposition}

\begin{Remark}
{\rm
After we recall the notion of local index in $\S$ 3, we will see that the condition $a_5=2a_1$ is equivalent to the following condition $(A)$:

$(A)$ \ The differential system $(2)$ has local index $(-a_1,-a_1,-a_1)$ at the accessible singular point $P_3=\{[z_0:z_1:z_2:z_3]=[0:1:0:0]\} \in {\cal H} \subset {\Bbb P}^3$. Here $[z_0:z_1:z_2:z_3]$ is the homogeneous coordinate of ${\Bbb P}^3$.
\rm}
\end{Remark}

\begin{Remark}
{\rm
In 1998, Noumi and Yamada [2] proposed the system of type $A_4^{(1)}$, which can be considered as a 4-parameter family of fourth-order coupled Painlev\'e IV system in dimension 4, and which is given as follows:

\begin{equation*}
  \left\{
  \begin{aligned}
   \frac{dx}{dt} &=x^2+2xy+2xw-tx+\beta_1\\
   \frac{dy}{dt} &=-y^2-2xy-2yw+ty+\beta_2\\
   \frac{dz}{dt} &=z^2+2zw+2xy-tz+\beta_1+\beta_3\\
   \frac{dw}{dt} &=-w^2-2zw+tw+\beta_4.
   \end{aligned}
  \right. 
\end{equation*}
This system reduces to the system $(1)$ by setting $x=0$ and $\beta_1=0$.
\rm}
\end{Remark}

In order to consider the phase space for the system $(1)$, let us take the compactification ${\Bbb P}^3$ of ${\Bbb C}^3$ with the natural embedding $(x,y,z)=(z_1/z_0,z_2/z_0,z_3/z_0)$. Moreover we denote the hypersuface of ${\Bbb P}^3$ by $ {\cal H}=\{z_0=0\} \simeq {\Bbb P}^2 \subset {\Bbb P}^3$. Fixing parameter $\alpha_i$ and setting $B={\rm Spec \rm} {\Bbb C} [t] \cong {\Bbb C}$, consider the product ${\Bbb P}^3 \times B \cong {\Bbb P}^3 \times \Bbb C$ and extend the regular vector field on ${\Bbb C}^3 \times B$ to a rational vector field $\tilde v$ on ${\Bbb P}^3 \times B$. By direct calculation, this rational vector field $\tilde v$ has seven accessible singular points on the boundary divisor ${\cal H} \times \{t\} \subset {\Bbb P}^3 \times \{t\}$ for each $t \in B$ (See Figure 1).

\begin{figure}[ht]
\unitlength 0.1in
\begin{picture}(32.60,15.35)(26.30,-25.65)
%
\special{pn 20}%
\special{pa 2828 1262}%
\special{pa 5413 1262}%
\special{fp}%
\special{pa 2837 1257}%
\special{pa 4179 1700}%
\special{fp}%
\special{pa 5428 1262}%
\special{pa 4161 1692}%
\special{fp}%
\special{pa 4171 1692}%
\special{pa 4171 2549}%
\special{fp}%
\special{pa 2846 1257}%
\special{pa 4154 2558}%
\special{fp}%
\special{pa 5438 1263}%
\special{pa 4179 2565}%
\special{fp}%
\put(45.4400,-26.1000){\makebox(0,0)[lb]{${\Bbb C}^3$}}%
\put(44.5800,-24.6200){\makebox(0,0)[lb]{$x$}}%
\put(42.0400,-24.1300){\makebox(0,0)[lb]{$y$}}%
\put(39.5000,-24.8300){\makebox(0,0)[lb]{$z$}}%
%
\special{pn 20}%
\special{sh 0.600}%
\special{ar 5438 1262 30 14  0.0000000 6.2831853}%
%
\special{pn 20}%
\special{sh 0.600}%
\special{ar 4212 1262 31 14  0.0000000 6.2831853}%
%
\special{pn 20}%
\special{sh 0.600}%
\special{ar 4169 1694 31 13  0.0000000 6.2831853}%
%
\special{pn 20}%
\special{sh 0.600}%
\special{ar 3493 1476 31 14  0.0000000 6.2831853}%
\put(39.8200,-17.2000){\makebox(0,0)[lb]{$P_6$}}%
\put(45.8700,-15.1100){\makebox(0,0)[lb]{$P_5$}}%
\put(52.6000,-12.0000){\makebox(0,0)[lb]{$P_3$}}%
\put(39.9000,-12.5000){\makebox(0,0)[lb]{$P_2$}}%
\put(26.3000,-12.2000){\makebox(0,0)[lb]{$P_1$}}%
\put(33.1200,-15.2200){\makebox(0,0)[lb]{$P_4$}}%
\put(39.9100,-14.5700){\makebox(0,0)[lb]{$P_7$}}%
\put(58.9000,-20.5000){\makebox(0,0)[lb]{${\Bbb P}^3$}}%
%
\special{pn 20}%
\special{pa 4160 2563}%
\special{pa 3960 2383}%
\special{fp}%
\special{sh 1}%
\special{pa 3960 2383}%
\special{pa 3996 2442}%
\special{pa 4000 2419}%
\special{pa 4023 2413}%
\special{pa 3960 2383}%
\special{fp}%
%
\special{pn 20}%
\special{pa 4170 2563}%
\special{pa 4170 2329}%
\special{fp}%
\special{sh 1}%
\special{pa 4170 2329}%
\special{pa 4150 2396}%
\special{pa 4170 2382}%
\special{pa 4190 2396}%
\special{pa 4170 2329}%
\special{fp}%
%
\special{pn 20}%
\special{pa 4170 2554}%
\special{pa 4380 2374}%
\special{fp}%
\special{sh 1}%
\special{pa 4380 2374}%
\special{pa 4316 2402}%
\special{pa 4340 2409}%
\special{pa 4342 2433}%
\special{pa 4380 2374}%
\special{fp}%
%
\special{pn 20}%
\special{ar 2850 1258 51 46  0.0000000 6.2831853}%
%
\special{pn 20}%
\special{ar 4590 1546 51 46  0.0000000 6.2831853}%
\put(38.8000,-15.3700){\makebox(0,0)[lb]{$\star$}}%
\end{picture}%
\label{fig:chazy2}
\caption{} 
\end{figure}

\begin{Theorem}

After a series of explicit blowing-ups and blowing-downs of ${\Bbb P}^3 \times B$, we obtain a smooth projective family of $3$-fold $\pi : {\mathcal X} \lra B$ and a birational morphism $\varphi : {\mathcal X} {\cdots \lra} {\Bbb P}^3 \times B$ which make the following diagram commutative$:$

$$
\begin{array}{ccc}
{\mathcal X} & \stackrel {\varphi} {\cdots \ra} & {\Bbb P}^3 \times B\\ 
\pi \downarrow \quad & & \downarrow\\ B & =& B,\\
\end{array}
$$
and satisfy the following conditions$:$

\begin{enumerate}

\item The phase space $\cal X$ over $B={\Bbb C}$ for the vector field $v$ in $(1)$ is obtained by gluing eight copies of ${\Bbb C}^3 \times {\Bbb C}:$

\begin{center}
${U_0} \times {\Bbb C}={{\Bbb C}^3} \times {\Bbb C} \ni {(x,y,z,t)},$

${U_j} \times {\Bbb C}={{\Bbb C}^3} \times {\Bbb C} \ni {(x_j,y_j,z_j,t)}$ $(j=1,2,....,7),$
\end{center}

\vspace{0.2cm}

\noindent
via the following rational transformations$:$

\vspace{0.2cm}

$1)$ $x_1=\frac{1}{x+z}, \ y_1=(yz+\alpha_2)z, \ z_1=\frac{1}{z},$

$2)$ $x_2=x+z,\ y_2=(yz+\alpha_2)z,\ z_2=\frac{1}{z}$,

$3)$ $x_3=\frac{1}{x}, \ y_3=y, \ z_3=z,$

$4)$ $x_4=(xz-\alpha_1)z, \ y_4=((y+z-t)z+1-\alpha_2-\alpha_3)z, \ z_4=\frac{1}{z},$

$5)$ $x_5=\frac{1}{x}, \ y_5=\frac{1}{y}, \ z_5=(yz-\alpha_3)y,$

$6)$ $x_6=x, \ y_6=\frac{1}{y}, \ z_6=(yz-\alpha_3)y,$

$7)$ $x_7=\frac{1}{z(xz-\alpha_1)},\ y_7=z((y+z-t)z+1-\alpha_2-\alpha_3), \ z_7=\frac{1}{z}$.

\item Each fiber $({\cX} \setminus {\cD})_t$ of the morphism $\pi : {\cX} 
\setminus {\cD} \lra B$ has a stratification

$$
({\cX} \setminus {\cD})_t= {\Bbb C}^3 \sqcup {\Bbb C}^2 \sqcup 
{\Bbb C}^2 \sqcup 
{\Bbb C}^2 \sqcup {\Bbb C}^2 \sqcup 
{\Bbb C}^1 \sqcup {\Bbb C}^1 \sqcup {\Bbb C}^1,
$$
where ${\Bbb C}^3$ is the original affine open subset of ${\Bbb P}^3 \setminus {\cal H}$. The three-dimensional space is a parameter space of the holomorphic solutions, and each two-dimensional or one-dimensional space is a parameter space of a two-parameter or a one-parameter family of meromorphic solutions, respectively.

\end{enumerate}
\end{Theorem}

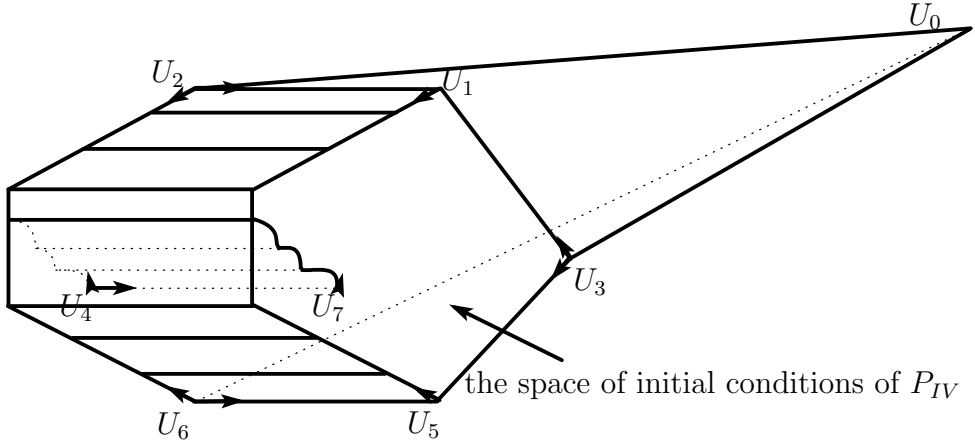
\begin{figure}[ht]
\unitlength 0.1in
\begin{picture}(50.40,21.51)(13.80,-24.52)
%
\special{pn 20}%
\special{pa 1384 1336}%
\special{pa 1384 1941}%
\special{fp}%
%
\special{pn 20}%
\special{pa 1384 1336}%
\special{pa 2659 1336}%
\special{fp}%
%
\special{pn 20}%
\special{pa 2659 1340}%
\special{pa 2659 1947}%
\special{fp}%
%
\special{pn 20}%
\special{pa 1389 1947}%
\special{pa 2664 1947}%
\special{fp}%
%
\special{pn 20}%
\special{pa 1384 1336}%
\special{pa 2364 808}%
\special{fp}%
%
\special{pn 20}%
\special{pa 2367 811}%
\special{pa 3644 811}%
\special{fp}%
%
\special{pn 20}%
\special{pa 2659 1340}%
\special{pa 3638 811}%
\special{fp}%
%
\special{pn 20}%
\special{pa 1776 1124}%
\special{pa 3050 1124}%
\special{fp}%
%
\special{pn 20}%
\special{pa 2139 934}%
\special{pa 3414 934}%
\special{fp}%
%
\special{pn 20}%
\special{pa 2664 1940}%
\special{pa 3634 2436}%
\special{fp}%
%
\special{pn 20}%
\special{pa 1380 1947}%
\special{pa 2351 2443}%
\special{fp}%
%
\special{pn 20}%
\special{pa 2354 2447}%
\special{pa 3629 2447}%
\special{fp}%
%
\special{pn 20}%
\special{pa 1721 2114}%
\special{pa 2996 2114}%
\special{fp}%
%
\special{pn 20}%
\special{pa 2081 2302}%
\special{pa 3356 2302}%
\special{fp}%
%
\special{pn 20}%
\special{pa 3648 808}%
\special{pa 4326 1700}%
\special{fp}%
\special{pa 3629 2443}%
\special{pa 4322 1700}%
\special{fp}%
%
\special{pn 20}%
\special{pa 1394 1496}%
\special{pa 2669 1496}%
\special{fp}%
%
\special{pn 8}%
\special{pa 1400 1496}%
\special{pa 1433 1505}%
\special{pa 1464 1516}%
\special{pa 1490 1532}%
\special{pa 1510 1555}%
\special{pa 1523 1583}%
\special{pa 1532 1616}%
\special{pa 1536 1633}%
\special{sp -0.045}%
%
\special{pn 8}%
\special{pa 1536 1642}%
\special{pa 1568 1653}%
\special{pa 1596 1668}%
\special{pa 1616 1690}%
\special{pa 1628 1719}%
\special{pa 1636 1752}%
\special{pa 1638 1762}%
\special{sp -0.045}%
%
\special{pn 8}%
\special{pa 1632 1759}%
\special{pa 1670 1759}%
\special{pa 1706 1761}%
\special{pa 1739 1766}%
\special{pa 1767 1777}%
\special{pa 1789 1795}%
\special{pa 1804 1821}%
\special{pa 1814 1852}%
\special{pa 1817 1865}%
\special{sp -0.045}%
%
\special{pn 8}%
\special{pa 1528 1646}%
\special{pa 2804 1646}%
\special{dt 0.045}%
\special{pa 2804 1646}%
\special{pa 2803 1646}%
\special{dt 0.045}%
%
\special{pn 8}%
\special{pa 1651 1759}%
\special{pa 2926 1759}%
\special{dt 0.045}%
\special{pa 2926 1759}%
\special{pa 2925 1759}%
\special{dt 0.045}%
%
\special{pn 8}%
\special{pa 1817 1854}%
\special{pa 3093 1854}%
\special{dt 0.045}%
\special{pa 3093 1854}%
\special{pa 3092 1854}%
\special{dt 0.045}%
%
\special{pn 20}%
\special{pa 2664 1491}%
\special{pa 2697 1504}%
\special{pa 2728 1518}%
\special{pa 2755 1535}%
\special{pa 2775 1557}%
\special{pa 2789 1584}%
\special{pa 2795 1615}%
\special{pa 2799 1646}%
\special{sp}%
%
\special{pn 20}%
\special{pa 2799 1642}%
\special{pa 2836 1639}%
\special{pa 2870 1642}%
\special{pa 2898 1659}%
\special{pa 2915 1696}%
\special{pa 2921 1747}%
\special{pa 2921 1759}%
\special{pa 2921 1754}%
\special{sp}%
%
\special{pn 20}%
\special{pa 2921 1759}%
\special{pa 2963 1756}%
\special{pa 3003 1754}%
\special{pa 3039 1756}%
\special{pa 3070 1762}%
\special{pa 3093 1776}%
\special{pa 3108 1797}%
\special{pa 3113 1827}%
\special{pa 3113 1859}%
\special{sp}%
%
\special{pn 20}%
\special{pa 2363 2447}%
\special{pa 2235 2383}%
\special{fp}%
\special{sh 1}%
\special{pa 2235 2383}%
\special{pa 2286 2431}%
\special{pa 2283 2407}%
\special{pa 2304 2395}%
\special{pa 2235 2383}%
\special{fp}%
%
\special{pn 20}%
\special{pa 2363 2447}%
\special{pa 2586 2444}%
\special{fp}%
\special{sh 1}%
\special{pa 2586 2444}%
\special{pa 2519 2425}%
\special{pa 2533 2445}%
\special{pa 2520 2465}%
\special{pa 2586 2444}%
\special{fp}%
%
\special{pn 20}%
\special{pa 3640 2441}%
\special{pa 3531 2389}%
\special{fp}%
\special{sh 1}%
\special{pa 3531 2389}%
\special{pa 3583 2436}%
\special{pa 3579 2412}%
\special{pa 3600 2400}%
\special{pa 3531 2389}%
\special{fp}%
%
\special{pn 20}%
\special{pa 4321 1693}%
\special{pa 4242 1788}%
\special{fp}%
\special{sh 1}%
\special{pa 4242 1788}%
\special{pa 4300 1750}%
\special{pa 4276 1747}%
\special{pa 4269 1724}%
\special{pa 4242 1788}%
\special{fp}%
%
\special{pn 20}%
\special{pa 4326 1693}%
\special{pa 4269 1605}%
\special{fp}%
\special{sh 1}%
\special{pa 4269 1605}%
\special{pa 4288 1672}%
\special{pa 4298 1650}%
\special{pa 4322 1650}%
\special{pa 4269 1605}%
\special{fp}%
%
\special{pn 20}%
\special{pa 2356 808}%
\special{pa 2235 872}%
\special{fp}%
\special{sh 1}%
\special{pa 2235 872}%
\special{pa 2303 859}%
\special{pa 2282 847}%
\special{pa 2285 823}%
\special{pa 2235 872}%
\special{fp}%
%
\special{pn 20}%
\special{pa 2363 804}%
\special{pa 2594 804}%
\special{fp}%
\special{sh 1}%
\special{pa 2594 804}%
\special{pa 2527 784}%
\special{pa 2541 804}%
\special{pa 2527 824}%
\special{pa 2594 804}%
\special{fp}%
%
\special{pn 20}%
\special{pa 3640 808}%
\special{pa 3518 872}%
\special{fp}%
\special{sh 1}%
\special{pa 3518 872}%
\special{pa 3586 859}%
\special{pa 3565 847}%
\special{pa 3568 823}%
\special{pa 3518 872}%
\special{fp}%
%
\special{pn 20}%
\special{pa 1817 1852}%
\special{pa 1805 1798}%
\special{fp}%
\special{sh 1}%
\special{pa 1805 1798}%
\special{pa 1800 1867}%
\special{pa 1817 1850}%
\special{pa 1839 1859}%
\special{pa 1805 1798}%
\special{fp}%
%
\special{pn 20}%
\special{pa 1825 1852}%
\special{pa 2023 1852}%
\special{fp}%
\special{sh 1}%
\special{pa 2023 1852}%
\special{pa 1956 1832}%
\special{pa 1970 1852}%
\special{pa 1956 1872}%
\special{pa 2023 1852}%
\special{fp}%
%
\special{pn 20}%
\special{pa 3113 1859}%
\special{pa 3126 1807}%
\special{fp}%
\special{sh 1}%
\special{pa 3126 1807}%
\special{pa 3090 1867}%
\special{pa 3113 1859}%
\special{pa 3129 1877}%
\special{pa 3126 1807}%
\special{fp}%
%
\special{pn 20}%
\special{pa 2375 808}%
\special{pa 6420 493}%
\special{fp}%
\put(60.9400,-4.7100){\makebox(0,0)[lb]{$U_0$}}%
%
\special{pn 20}%
\special{pa 6414 498}%
\special{pa 4326 1697}%
\special{fp}%
%
\special{pn 8}%
\special{pa 6414 494}%
\special{pa 2362 2449}%
\special{dt 0.045}%
\special{pa 2362 2449}%
\special{pa 2363 2449}%
\special{dt 0.045}%
\put(21.3000,-7.7300){\makebox(0,0)[lb]{$U_2$}}%
\put(36.5000,-8.1300){\makebox(0,0)[lb]{$U_1$}}%
\put(43.4000,-18.6200){\makebox(0,0)[lb]{$U_3$}}%
\put(34.7000,-26.1000){\makebox(0,0)[lb]{$U_5$}}%
\put(21.6000,-26.1800){\makebox(0,0)[lb]{$U_6$}}%
\put(16.5000,-20.0000){\makebox(0,0)[lb]{$U_4$}}%
\put(29.7000,-20.0000){\makebox(0,0)[lb]{$U_7$}}%
%
\special{pn 20}%
\special{pa 4280 2230}%
\special{pa 3720 1954}%
\special{fp}%
\special{sh 1}%
\special{pa 3720 1954}%
\special{pa 3771 2001}%
\special{pa 3768 1978}%
\special{pa 3789 1966}%
\special{pa 3720 1954}%
\special{fp}%
\put(37.7000,-24.0100){\makebox(0,0)[lb]{the space of initial conditions of $P_{IV}$}}%
\end{picture}%
\label{fig:chazy17}
\caption{}
\end{figure}

\begin{Remark}
{\rm
The surface defined by $x=0$ is the minimal compactification of the space of initial conditions of Painlev\'e IV equation (see Figure 2). Actually, the ordinary differential system $(1)$ can be reduced to Painlev\'e IV equation if one sets $x \equiv 0$ and $\alpha_1=0$. However, the total differential system (1) is not a product of the fourth Painlev\'e equation and a Riccati equation.

\rm}
\end{Remark}

\begin{Remark}
{\rm
It is still an open question whether we can resolve the accessible singularities for any differential system in $(2)$ with the condition $a_5=2a_1$.
\rm}
\end{Remark}

\section{ Accessible singularity}

\vspace{0.2cm}

Let us review the notion of accessible singularity [4]. Let $B$ be a connected open domain in $\Bbb C$ and $\pi : {\mathcal W} \lra B$ a smooth proper holomorphic map. We assume that ${\mathcal H} \subset {\mathcal W}$ is a normal crossing divisor which is flat over $B$. Let us consider a rational vector field $\tilde v$ on $\mathcal W$ satisfying the condition

$$
\tilde v \in H^0({\cal W},\Theta_{\cal W}(-\log{\cH})({\cH})).
$$

\noindent
Fixing $t_0 \in B$ and $P \in {\mathcal W}_{t_0}$, we can take a local coordinate system $(x_1,x_2,...,x_n)$ of ${\cal W}_{t_0}$ centered at $P$ such that ${\cH}_{smooth}$ can be defined by the local equation $x_1=0$.
Since $\tilde v \in H^0({\cal W},\Theta_{\cal W}(-\log{\cH})({\cH}))$, we can write down the vector field $\tilde v$ near $P=(0,0,...,0,t_0)$ as follows:

\begin{equation*}
\boxed{%
\tilde v= \frac{\partial}{\partial t}+a_1 
\frac{\partial}{\partial x_1}+\frac{a_2}{x_1} 
\frac{\partial}{\partial x_2}+.....+\frac{a_n}{x_1} 
\frac{\partial}{\partial x_n}
}%
\end{equation*}

\noindent
This vector field defines the following system of differential equations

\begin{equation}
  \left\{
  \begin{aligned}
   \frac{dx_1}{dt} &=a_1(x_1,x_2,....,x_n,t)\\
   \frac{dx_2}{dt} &=\frac{a_2(x_1,x_2,....,x_n,t)}{x_1}\\
   .\\
   .\\
   .\\
   \frac{dx_n}{dt} &=\frac{a_n(x_1,x_2,....,x_n,t)}{x_1}.
   \end{aligned}
  \right. 
\end{equation}

\noindent
Here $a_i(x_1,...,x_n,t)$ $(i=1,2,...,n)$ are holomorphic functions near $P=(0,0,..,0,t_0).$

\begin{Definition}

With the above notation, assume that the rational vector field $\tilde v$ on $\cal W$ satisfies the condition
$$
\tilde v \in H^0({\cal W},\Theta_{\cal W}(-\log{\cH})({\cH})).
$$
We say that $\tilde v$ has an accessible singular point at $P=(0,0,...,0,t_0)$ if

$$
x_1=0 \ {\rm and \rm} \ a_i(0,0,....,0,t_0)=0 \ {\rm for \rm} \ {\rm every \rm} \ i, \ 2 \leq i \leq n.
$$

\end{Definition}

If $P \in {\cH}_{smooth}$ is not an accessible singular point, all solutions of the ordinary differential equation passing through $P$ are vertical solutions, that is, the solutions are contained in the fiber ${\cal W}_{t_0}$ over $t=t_0$. If $P \in {\cH}_{smooth}$ is an accessible singular point, there may be a solution of $(3)$ which passes through $P$ and goes into the interior ${\cal W}-{\cH}$ of ${\cal W}$.

Let us recall the notion of local index [4]. When we construct the phase spaces of the higher order Painlev\'e equations, an object that we call the local index is the key for determining when we need to make a blowing-up of an accessible singularity or a blowing-down to a minimal phase space. In the case of equations of higher order with favorable properties, for example the systems of type ${A_4}^{(1)}$ [2], the local index at the accessible singular point corresponds to the set of orders that appears in the free parameters of formal solutions passing through that point [5].

\vspace{0.2cm}

\begin{Definition}
Let $v$ be an algebraic vector field which is given by $(3)$ and $(X,Y,Z)$ be a boundary coordinate system in a neighborhood of an accessible singular point $P=(0,0,0,t)$. Assume that the system is written as

\begin{equation*}
  \left\{
\begin{aligned}
   \frac {dX}{dt} &=a+f_1(X,Y,Z,t)\\
   \frac {dY}{dt} &=\frac{bY+f_2(X,Y,Z,t)}{X}\\
   \frac {dZ}{dt} &=\frac{cZ+f_3(X,Y,Z,t)}{X}
\end{aligned}
  \right. 
\end{equation*}
near the accessible singular point $P$, where $a,b$ and $c$ are nonzero constants. We say that the vector field $v$ has local index $(a,b,c)$ at $P$ if $f_1(X,Y,Z,t)$ is a polynomial which vanishes at $P=(0,0,0,t)$ and $f_i(X,Y,Z,t)$ $(i=2,3)$ are polynomials of order $2$ in $X,Y,Z$. Here $f_i \in {\Bbb C}[X,Y,Z,t]$ $(i=1,2,3).$

\end{Definition}

\vspace{0.2cm}
\begin{Remark}
{\rm
We are interested in the case with local index $(1,\frac{b}{a},\frac{c}{a}) \in {\Bbb Z}^3$. If each component of $(1,\frac{b}{a},\frac{c}{a})$ has the same sign, we may resolve the accessible singularity by blowing-up finitely many times. But when different signs appear, we may need to both blow up and blow down.
\rm}
\end{Remark}

In our case, there exist 7 accessible singular points of this system on the boundary divisor of ${\Bbb P}^3$.  They are listed as follows.

\vspace{0.5cm}

\begin{center}
\begin{tabular}{|c||c|c|c|c|} \hline 
Singular point & $P_1$ & $P_2$ & $P_3$ & $P_4$ \\ \hline 
$[z_0:z_1:z_2:z_3]$& $[0: 0: 0: 1]$& $[0: -1: 0: 1] $& $[0:1:0:0]$ & $[0:0: -1:1]$ \\ \hline 
Type (dim. of sol.)    &  $\circ$ (dim. 1) & $\bullet$ (dim. 2)  & $\bullet$ (dim. 2)  &
$\bullet$ (dim. 2)  \\ \hline
Type of local index & (-1,+3,+1) & (+1,+3,+1) & (+1,+1,+1) & (-3,-3,-1)   \\ \hline
\end{tabular}

\begin{tabular}{|c||c|c|c|} \hline 
Singular point  & $P_5$  & $P_6$ &  $P_7$ \\ \hline 
$[z_0:z_1:z_2:z_3]$   &$[0:-1: 1:0]$ & $[0:0:1:0]$ & $[0:-3:-1:1]$ \\  \hline 
Type (dim. of sol.)    & $\circ$ (dim. 1)    &  $\bullet$ (dim. 2)     & $\star$ (dim. 1)  \\ \hline
Type of local index & (+1,-1,-3) & (-1,-1,-3) & (+3,-3,-1)      \\ \hline
\end{tabular}
\end{center}

\vspace{0.5cm}

\noindent
Here, we remark that there are 3 types of accessible singularities.  
We denote the mark by $\bullet$ an accessible singular point into which two dimensional meromorphic solutions flow, and the marks by $\circ$ and $\star$ accessible singular points into which one dimensional meromorphic solutions flow.

\section{ Resolution of accessible singularities}

Let $P$ be an accessible singular point in the boundary divisor $\cal H$. Rewriting the systems in a local coordinate system, the right hand side of each differential equation has poles along $\cal H$. If we resolve the accessible singular point $P$ and the right hand side of each differential equation becomes holomorphic in new coordinate system, then we can use Cauchy's existence and uniqueness theorem of solutions.

{\bf Notation}

$(u,v,w)=(\frac{x}{z},\frac{y}{z},\frac{1}{z}), \ (p,q,r)=(\frac{1}{x},\frac{y}{x},\frac{z}{x}), \ (l,m,n)=(\frac{x}{y},\frac{1}{y},\frac{z}{y}), \ (\tilde u,v,w)=(u+1,v,w).$

{\bf 4.1. Resolution of accessible singular points $P_1,P_2$ and $P_3$ }

Let us start by summarizing the steps which are needed to resolve the accessible singular points of $\tilde v$ at $ P_1:(u,v,w)=(0,0,0),  P_2:(\tilde u,v,w)=(0,0,0) \ and \ P_3:(p,q,r)=(0,0,0)$.

\begin{enumerate}
\item{\bf Step 1}: We will blow up at two points $P_2, P_3$ (see Figure 3).
\item{\bf Step 2}: We will blow up along the curve $C_1$ (see Figure 3).
\item{\bf Step 3}: We will blow down the surface $F$ (see Figure 3).
\item{\bf Step 4}: We will blow up along the curve $C_4$ (see Figure 3).
\item{\bf Step 5}: We will blow up along the curve $C_5$ (see Figure 3).
\end{enumerate}

\begin{figure}[ht]
\unitlength 0.1in
\begin{picture}(54.72,25.17)(8.60,-25.97)
%
\special{pn 20}%
\special{pa 4938 82}%
\special{pa 4343 248}%
\special{fp}%
\special{pa 4343 248}%
\special{pa 4764 394}%
\special{fp}%
\special{pa 4938 80}%
\special{pa 5612 247}%
\special{fp}%
\special{pa 5612 247}%
\special{pa 5198 394}%
\special{fp}%
\special{pa 4773 394}%
\special{pa 5202 394}%
\special{fp}%
\special{pa 4348 248}%
\special{pa 4179 491}%
\special{fp}%
\special{pa 4179 491}%
\special{pa 4769 394}%
\special{fp}%
\special{pa 5612 247}%
\special{pa 5926 489}%
\special{fp}%
\special{pa 5202 392}%
\special{pa 5930 489}%
\special{fp}%
%
\special{pn 20}%
\special{pa 4184 487}%
\special{pa 4340 604}%
\special{fp}%
\special{sh 1}%
\special{pa 4340 604}%
\special{pa 4299 548}%
\special{pa 4297 572}%
\special{pa 4275 580}%
\special{pa 4340 604}%
\special{fp}%
\special{pa 5935 489}%
\special{pa 5712 628}%
\special{fp}%
\special{sh 1}%
\special{pa 5712 628}%
\special{pa 5779 610}%
\special{pa 5757 600}%
\special{pa 5758 576}%
\special{pa 5712 628}%
\special{fp}%
%
\special{pn 20}%
\special{sh 0.600}%
\special{ar 4188 490 49 15  0.0000000 6.2831853}%
%
\special{pn 20}%
\special{pa 2218 158}%
\special{pa 1520 540}%
\special{fp}%
\special{pa 2212 158}%
\special{pa 2921 542}%
\special{fp}%
%
\special{pn 20}%
\special{pa 1514 536}%
\special{pa 2909 536}%
\special{fp}%
%
\special{pn 20}%
\special{sh 0.600}%
\special{ar 1526 540 48 13  0.0000000 6.2831853}%
%
\special{pn 20}%
\special{sh 0.600}%
\special{ar 2909 540 49 13  0.0000000 6.2831853}%
\put(19.1000,-5.2000){\makebox(0,0)[lb]{$P_1$}}%
\put(12.0600,-6.1900){\makebox(0,0)[lb]{$P_2$}}%
\put(28.2000,-4.5000){\makebox(0,0)[lb]{$P_3$}}%
\put(49.4000,-7.0000){\makebox(0,0)[lb]{$C_1$}}%
\put(34.0800,-3.1100){\makebox(0,0){Step 1}}%
%
\special{pn 20}%
\special{pa 3646 389}%
\special{pa 3102 389}%
\special{fp}%
\special{sh 1}%
\special{pa 3102 389}%
\special{pa 3169 409}%
\special{pa 3155 389}%
\special{pa 3169 369}%
\special{pa 3102 389}%
\special{fp}%
%
\special{pn 20}%
\special{pa 1325 1193}%
\special{pa 900 1193}%
\special{fp}%
\special{sh 1}%
\special{pa 900 1193}%
\special{pa 967 1213}%
\special{pa 953 1193}%
\special{pa 967 1173}%
\special{pa 900 1193}%
\special{fp}%
%
\special{pn 20}%
\special{pa 2328 872}%
\special{pa 1559 1101}%
\special{fp}%
\special{pa 1559 1106}%
\special{pa 1486 1423}%
\special{fp}%
\special{pa 1559 1101}%
\special{pa 1997 1225}%
\special{fp}%
\special{pa 2004 1225}%
\special{pa 2004 1381}%
\special{fp}%
%
\special{pn 20}%
\special{pa 1493 1423}%
\special{pa 2019 1376}%
\special{fp}%
%
\special{pn 20}%
\special{pa 2321 866}%
\special{pa 3019 1101}%
\special{fp}%
\special{pa 3019 1094}%
\special{pa 2558 1235}%
\special{fp}%
\special{pa 2558 1235}%
\special{pa 2558 1365}%
\special{fp}%
\special{pa 2558 1365}%
\special{pa 3120 1418}%
\special{fp}%
\special{pa 3012 1094}%
\special{pa 3127 1423}%
\special{fp}%
%
\special{pn 20}%
\special{pa 2004 1225}%
\special{pa 2558 1225}%
\special{fp}%
%
\special{pn 20}%
\special{pa 2012 1376}%
\special{pa 2558 1376}%
\special{fp}%
%
\special{pn 20}%
\special{pa 3322 1216}%
\special{pa 3875 1216}%
\special{fp}%
\special{sh 1}%
\special{pa 3875 1216}%
\special{pa 3808 1196}%
\special{pa 3822 1216}%
\special{pa 3808 1236}%
\special{pa 3875 1216}%
\special{fp}%
%
\special{pn 20}%
\special{pa 4734 880}%
\special{pa 4358 1071}%
\special{fp}%
\special{pa 4358 1071}%
\special{pa 5311 1160}%
\special{fp}%
\special{pa 4740 872}%
\special{pa 5699 964}%
\special{fp}%
\special{pa 5706 964}%
\special{pa 5298 1155}%
\special{fp}%
%
\special{pn 20}%
\special{pa 4364 1071}%
\special{pa 4232 1479}%
\special{fp}%
\special{pa 4232 1479}%
\special{pa 5298 1403}%
\special{fp}%
\special{pa 5304 1160}%
\special{pa 5304 1387}%
\special{fp}%
%
\special{pn 20}%
\special{pa 5712 960}%
\special{pa 6144 1176}%
\special{fp}%
\special{pa 5317 1387}%
\special{pa 6144 1176}%
\special{fp}%
%
\special{pn 20}%
\special{pa 1342 2069}%
\special{pa 917 2069}%
\special{fp}%
\special{sh 1}%
\special{pa 917 2069}%
\special{pa 984 2089}%
\special{pa 970 2069}%
\special{pa 984 2049}%
\special{pa 917 2069}%
\special{fp}%
%
\special{pn 20}%
\special{pa 1903 1774}%
\special{pa 1527 1965}%
\special{fp}%
\special{pa 1527 1965}%
\special{pa 2480 2053}%
\special{fp}%
\special{pa 1910 1766}%
\special{pa 2869 1857}%
\special{fp}%
\special{pa 2875 1857}%
\special{pa 2467 2049}%
\special{fp}%
%
\special{pn 20}%
\special{pa 1533 1965}%
\special{pa 1402 2372}%
\special{fp}%
\special{pa 1402 2372}%
\special{pa 2467 2296}%
\special{fp}%
\special{pa 2473 2053}%
\special{pa 2473 2281}%
\special{fp}%
%
\special{pn 20}%
\special{pa 2881 1853}%
\special{pa 3314 2069}%
\special{fp}%
\special{pa 2486 2281}%
\special{pa 3314 2069}%
\special{fp}%
%
\special{pn 20}%
\special{pa 1478 2164}%
\special{pa 2481 2164}%
\special{fp}%
%
\special{pn 8}%
\special{pa 1452 2242}%
\special{pa 2472 2242}%
\special{dt 0.045}%
\special{pa 2472 2242}%
\special{pa 2471 2242}%
\special{dt 0.045}%
%
\special{pn 20}%
\special{pa 3884 2069}%
\special{pa 3458 2069}%
\special{fp}%
\special{sh 1}%
\special{pa 3458 2069}%
\special{pa 3525 2089}%
\special{pa 3511 2069}%
\special{pa 3525 2049}%
\special{pa 3458 2069}%
\special{fp}%
%
\special{pn 8}%
\special{pa 4496 1799}%
\special{pa 4172 1985}%
\special{dt 0.045}%
\special{pa 4172 1985}%
\special{pa 4173 1985}%
\special{dt 0.045}%
\special{pa 4172 1985}%
\special{pa 5260 2052}%
\special{dt 0.045}%
\special{pa 5260 2052}%
\special{pa 5259 2052}%
\special{dt 0.045}%
\special{pa 4496 1788}%
\special{pa 5754 1844}%
\special{dt 0.045}%
\special{pa 5754 1844}%
\special{pa 5753 1844}%
\special{dt 0.045}%
\special{pa 5754 1844}%
\special{pa 5278 2052}%
\special{dt 0.045}%
\special{pa 5278 2052}%
\special{pa 5279 2052}%
\special{dt 0.045}%
%
\special{pn 8}%
\special{pa 4181 1990}%
\special{pa 4088 2367}%
\special{dt 0.045}%
\special{pa 4088 2367}%
\special{pa 4088 2366}%
\special{dt 0.045}%
%
\special{pn 20}%
\special{pa 5745 1851}%
\special{pa 6332 2400}%
\special{fp}%
%
\special{pn 8}%
\special{pa 5269 2058}%
\special{pa 5269 2254}%
\special{dt 0.045}%
\special{pa 5269 2254}%
\special{pa 5269 2253}%
\special{dt 0.045}%
%
\special{pn 8}%
\special{pa 4172 2069}%
\special{pa 5269 2136}%
\special{dt 0.045}%
\special{pa 5269 2136}%
\special{pa 5268 2136}%
\special{dt 0.045}%
%
\special{pn 8}%
\special{pa 4130 2170}%
\special{pa 5269 2254}%
\special{dt 0.045}%
\special{pa 5269 2254}%
\special{pa 5268 2254}%
\special{dt 0.045}%
%
\special{pn 20}%
\special{pa 4138 2175}%
\special{pa 4170 2185}%
\special{pa 4203 2195}%
\special{pa 4234 2206}%
\special{pa 4264 2219}%
\special{pa 4293 2233}%
\special{pa 4320 2249}%
\special{pa 4344 2267}%
\special{pa 4366 2289}%
\special{pa 4385 2313}%
\special{pa 4403 2340}%
\special{pa 4418 2368}%
\special{pa 4432 2398}%
\special{pa 4444 2423}%
\special{sp}%
%
\special{pn 20}%
\special{pa 5260 2254}%
\special{pa 5292 2264}%
\special{pa 5325 2274}%
\special{pa 5356 2285}%
\special{pa 5386 2298}%
\special{pa 5415 2312}%
\special{pa 5442 2328}%
\special{pa 5466 2346}%
\special{pa 5488 2368}%
\special{pa 5507 2392}%
\special{pa 5525 2419}%
\special{pa 5540 2447}%
\special{pa 5554 2477}%
\special{pa 5566 2502}%
\special{sp}%
%
\special{pn 20}%
\special{pa 4436 2412}%
\special{pa 5566 2502}%
\special{fp}%
\special{pa 5566 2508}%
\special{pa 6323 2406}%
\special{fp}%
%
\special{pn 8}%
\special{pa 4079 2367}%
\special{pa 4402 2367}%
\special{dt 0.045}%
\special{pa 4402 2367}%
\special{pa 4401 2367}%
\special{dt 0.045}%
%
\special{pn 8}%
\special{pa 4419 2372}%
\special{pa 5269 2372}%
\special{dt 0.045}%
\special{pa 5269 2372}%
\special{pa 5268 2372}%
\special{dt 0.045}%
%
\special{pn 8}%
\special{pa 5269 2249}%
\special{pa 5269 2372}%
\special{dt 0.045}%
\special{pa 5269 2372}%
\special{pa 5269 2371}%
\special{dt 0.045}%
\put(11.3000,-11.0800){\makebox(0,0){Step 2}}%
\put(35.9400,-11.3100){\makebox(0,0){Step 3}}%
\put(11.3800,-19.9600){\makebox(0,0){Step 4}}%
\put(37.0500,-19.9000){\makebox(0,0){Step 5}}%
%
\special{pn 20}%
\special{ar 2252 1389 29 21  0.0000000 6.2831853}%
%
\special{pn 20}%
\special{ar 5312 1407 29 21  0.0000000 6.2831853}%
%
\special{pn 20}%
\special{ar 2464 2242 30 22  0.0000000 6.2831853}%
%
\special{pn 20}%
\special{pa 2252 1367}%
\special{pa 2252 1266}%
\special{fp}%
\special{sh 1}%
\special{pa 2252 1266}%
\special{pa 2232 1333}%
\special{pa 2252 1319}%
\special{pa 2272 1333}%
\special{pa 2252 1266}%
\special{fp}%
%
\special{pn 20}%
\special{pa 2252 1372}%
\special{pa 2472 1372}%
\special{fp}%
\special{sh 1}%
\special{pa 2472 1372}%
\special{pa 2405 1352}%
\special{pa 2419 1372}%
\special{pa 2405 1392}%
\special{pa 2472 1372}%
\special{fp}%
%
\special{pn 20}%
\special{pa 1520 536}%
\special{pa 1674 687}%
\special{fp}%
\special{sh 1}%
\special{pa 1674 687}%
\special{pa 1640 626}%
\special{pa 1636 650}%
\special{pa 1612 655}%
\special{pa 1674 687}%
\special{fp}%
\special{pa 2914 536}%
\special{pa 2778 692}%
\special{fp}%
\special{sh 1}%
\special{pa 2778 692}%
\special{pa 2837 655}%
\special{pa 2813 652}%
\special{pa 2807 629}%
\special{pa 2778 692}%
\special{fp}%
%
\special{pn 20}%
\special{pa 1486 1416}%
\special{pa 1640 1597}%
\special{fp}%
\special{sh 1}%
\special{pa 1640 1597}%
\special{pa 1612 1533}%
\special{pa 1605 1556}%
\special{pa 1582 1559}%
\special{pa 1640 1597}%
\special{fp}%
%
\special{pn 20}%
\special{pa 3118 1412}%
\special{pa 2906 1613}%
\special{fp}%
\special{sh 1}%
\special{pa 2906 1613}%
\special{pa 2968 1582}%
\special{pa 2945 1576}%
\special{pa 2941 1553}%
\special{pa 2906 1613}%
\special{fp}%
%
\special{pn 20}%
\special{pa 4240 1473}%
\special{pa 4470 1619}%
\special{fp}%
\special{sh 1}%
\special{pa 4470 1619}%
\special{pa 4424 1566}%
\special{pa 4425 1590}%
\special{pa 4403 1600}%
\special{pa 4470 1619}%
\special{fp}%
%
\special{pn 20}%
\special{pa 6136 1170}%
\special{pa 6034 1479}%
\special{fp}%
\special{sh 1}%
\special{pa 6034 1479}%
\special{pa 6074 1422}%
\special{pa 6051 1428}%
\special{pa 6036 1409}%
\special{pa 6034 1479}%
\special{fp}%
%
\special{pn 20}%
\special{pa 1418 2378}%
\special{pa 1708 2563}%
\special{fp}%
\special{sh 1}%
\special{pa 1708 2563}%
\special{pa 1663 2510}%
\special{pa 1663 2534}%
\special{pa 1641 2544}%
\special{pa 1708 2563}%
\special{fp}%
%
\special{pn 20}%
\special{pa 3297 2069}%
\special{pa 3212 2384}%
\special{fp}%
\special{sh 1}%
\special{pa 3212 2384}%
\special{pa 3249 2325}%
\special{pa 3226 2333}%
\special{pa 3210 2314}%
\special{pa 3212 2384}%
\special{fp}%
%
\special{pn 20}%
\special{pa 4453 2423}%
\special{pa 4691 2597}%
\special{fp}%
\special{sh 1}%
\special{pa 4691 2597}%
\special{pa 4649 2542}%
\special{pa 4648 2566}%
\special{pa 4625 2574}%
\special{pa 4691 2597}%
\special{fp}%
%
\special{pn 20}%
\special{pa 6332 2395}%
\special{pa 6298 2592}%
\special{fp}%
\special{sh 1}%
\special{pa 6298 2592}%
\special{pa 6329 2530}%
\special{pa 6307 2539}%
\special{pa 6290 2523}%
\special{pa 6298 2592}%
\special{fp}%
\put(22.6000,-14.7400){\makebox(0,0)[lb]{$u_1$}}%
\put(20.5000,-13.2400){\makebox(0,0)[lb]{$v_1$}}%
%
\special{pn 20}%
\special{pa 5300 1386}%
\special{pa 5610 1318}%
\special{fp}%
\special{sh 1}%
\special{pa 5610 1318}%
\special{pa 5541 1313}%
\special{pa 5558 1329}%
\special{pa 5549 1352}%
\special{pa 5610 1318}%
\special{fp}%
%
\special{pn 20}%
\special{pa 5280 1401}%
\special{pa 4810 1432}%
\special{fp}%
\special{sh 1}%
\special{pa 4810 1432}%
\special{pa 4878 1448}%
\special{pa 4863 1428}%
\special{pa 4875 1408}%
\special{pa 4810 1432}%
\special{fp}%
%
\special{pn 20}%
\special{pa 5280 1396}%
\special{pa 5280 1235}%
\special{fp}%
\special{sh 1}%
\special{pa 5280 1235}%
\special{pa 5260 1302}%
\special{pa 5280 1288}%
\special{pa 5300 1302}%
\special{pa 5280 1235}%
\special{fp}%
\put(51.1000,-12.8000){\makebox(0,0)[lb]{$v_2$}}%
\put(46.6000,-14.1600){\makebox(0,0)[lb]{$u_2$}}%
\put(53.8000,-13.6500){\makebox(0,0)[lb]{$w_2$}}%
%
\special{pn 20}%
\special{pa 5110 540}%
\special{pa 5120 390}%
\special{fp}%
\special{sh 1}%
\special{pa 5120 390}%
\special{pa 5096 455}%
\special{pa 5116 443}%
\special{pa 5136 458}%
\special{pa 5120 390}%
\special{fp}%
%
\special{pn 20}%
\special{ar 2140 560 45 28  0.0000000 6.2831853}%
%
\special{pn 20}%
\special{ar 4970 400 45 28  0.0000000 6.2831853}%
\end{picture}%
\label{fig:third1}
\caption{}
\end{figure}

In a neighborhood of $P_1$, the system is written as

\begin{equation*}
  \left\{
  \begin{aligned}
   \frac {du}{dt} &=\frac{-u-u^2+2uv}{w}+\alpha_1w-\alpha_3uw\\
   \frac {dv}{dt} &=-2tv+\frac{3v+3v^2}{w}+\alpha_2w-\alpha_3vw\\
   \frac {dw}{dt} &=1+2v-tw-\alpha_3w^2.
   \end{aligned}
  \right. 
\end{equation*}

In a neighborhood of $P_2$, the system is written as

\begin{equation*}
  \left\{
  \begin{aligned}
   \frac {d{\tilde u}}{dt} &=\frac{\tilde u-{\tilde u}^2-v(2-2\tilde u)}{w}+\alpha_3w+
   \alpha_1w-\alpha_3{\tilde u}w\\
   \frac {dv}{dt} &=-2tv+\frac{3v+3{v}^2}{w}+\alpha_2w-\alpha_3vw\\
   \frac {dw}{dt} &=-\alpha_3w^2-tw+2v+1.\\ 
   \end{aligned}
  \right. 
\end{equation*}

In a neighborhood of $P_3$, the system is written as

\begin{equation*}
  \left\{
  \begin{aligned}
   \frac {dp}{dt} &=-\alpha_1p^2-tp+2r+1\\
   \frac {dq}{dt} &=\frac{q+q(-2tp+q+4r)}{p}+\alpha_2p-\alpha_1pq\\
   \frac {dr}{dt} &=\frac{r+r(r-2q)}{p}-\alpha_1pr+\alpha_3p.\\
   \end{aligned}
  \right. 
\end{equation*}

{\bf 4.1.1.} We blow up at two points $P_2, P_3$.

We blow up at the point $P_2$

$$
U_1=\frac{\tilde u}{w} \;, \;\;\; V_1=\frac{v}{w} \;, \;\;\; W_1=w.
$$

We blow up at the point $P_3$ 

$$
x_3=p \;, \;\;\; y_3=\frac{q}{p} \;, \;\;\; z_3=\frac{r}{p}.
$$
We have resolved the accessible singular point $P_3$.

\vspace{0.2cm}

{\bf 4.1.2.} We blow up along the curve $C_1=\{(u,v,w)| v=w=0\} \subset {\Bbb P}^1$

$$
u_1=u \;, \;\;\; v_1=\frac{v}{w} \;, \;\;\; w_1=w.
$$
In a neighborhood of $\{(u_1,v_1,w_1)=(0,0,0)\}$, the system is written as

\begin{equation*}
  \left\{
  \begin{aligned}
   \frac {du_1}{dt} &=2u_1v_1+\frac{-u_1-(u_1)^2}{w_1}+\alpha_1w_1-\alpha_3
   u_1w_1\\
   \frac {dv_1}{dt} &=-tv_1+(v_1)^2+\frac{2v_1}{w_1}+\alpha_2\\
   \frac {dw_1}{dt} &=1-tw_1+2v_1w_1-\alpha_3(w_1)^2.
   \end{aligned}
  \right. 
\end{equation*}

{\bf 4.1.3.} We blow down the surface $F=\{(u_1,v_1,w_1)|w_1=0\} \cong {\Bbb P}^1 \times {\Bbb P}^1$

$$
u_2=\frac{w_1}{u_1+1} \;, \;\;\; v_2=v_1 \;, \;\;\; w_2=w_1.
$$
In a neighborhood of $\{(u_2,v_2,w_2)=(0,0,0)\}$, the system is written as

\begin{equation}
  \left\{
  \begin{aligned}
   \frac {du_2}{dt} &=1-tu_2+\frac{2(u_2)^2v_2}{w_2}-\alpha_1(u_2)^2
   -\alpha_3(u_2)^2\\
   \frac {dv_2}{dt} &=-tv_2+(v_2)^2+\frac{2v_2}{w_2}+\alpha_2\\
   \frac {dw_2}{dt} &=1-tw_2+2v_2w_2-\alpha_3(w_2)^2.
   \end{aligned}
  \right. 
\end{equation}

The resolution process from Step 2 to Step 3 is well-known as flop. In order to resolve the accessible singular point $P_1$ and obtain a holomorphic coordinate system, we need to blow down the surface $F \cong {\Bbb P}^1 \times {\Bbb P}^1$ along the ${\Bbb P}^1$-fiber. After we blow down the surface $F$, the system $(4)$ has the local index $(0,+2,+1)$ at the point $\{(u_2,v_2,w_2)=(0,0,0)\}$.

{\bf 4.1.4.} We blow up along the curve $C_4=\{(u_2,v_2,w_2)| v_2=w_2=0\} \subset {\Bbb P}^1$

$$
u_3=u_2 \;, \;\;\; v_3=\frac{v_2}{w_2} \;, \;\;\; w_3=w_2.
$$

{\bf 4.1.5.} We blow up along the curve $C_5=\{(u_3,v_3,w_3)| v_3+\alpha_2=w_3=0\} \subset {\Bbb P}^1$

$$
x_1=u_3 \;, \;\;\; y_1=\frac{v_3+\alpha_2}{w_3} \;, \;\;\; z_1=w_3.
$$
We have resolved the accessible singular point $P_1$.

\vspace{0.3cm}

{\bf 4.2. Resolution of accessible singular points $ P_5:(l,m,n)=(0,-1,0)$ and $P_6:(l,m,n)=(0,0,0)$ }

We can resolve accessible singular points $P_5, P_6$ by the same way of 4.1. By the birational transformation of 4.1 and 4.2, we have resolved five accessible singular points $P_1, P_2, P_3, P_5$ and $P_6$. There are only two accessible singular points $P_4, P_7$ in ${\Bbb F}_1$-surface. Here ${\Bbb F}_1$-surface is obtained by blowing up one time at a point in ${\Bbb P}^2$. These two accessible singular points are on the same fiber in ${\Bbb F}_1$-surface (See figure 4).

\begin{figure}[ht]
\unitlength 0.1in
\begin{picture}(57.72,41.10)(9.40,-45.10)
%
\special{pn 20}%
\special{pa 1688 736}%
\special{pa 1688 1120}%
\special{fp}%
%
\special{pn 20}%
\special{pa 1688 736}%
\special{pa 2617 736}%
\special{fp}%
%
\special{pn 20}%
\special{pa 2617 738}%
\special{pa 2617 1123}%
\special{fp}%
%
\special{pn 20}%
\special{pa 1691 1123}%
\special{pa 2621 1123}%
\special{fp}%
%
\special{pn 20}%
\special{pa 1688 736}%
\special{pa 2402 400}%
\special{fp}%
%
\special{pn 20}%
\special{pa 2405 403}%
\special{pa 3333 403}%
\special{fp}%
%
\special{pn 20}%
\special{pa 2617 738}%
\special{pa 3331 403}%
\special{fp}%
%
\special{pn 20}%
\special{pa 1971 601}%
\special{pa 2902 601}%
\special{fp}%
%
\special{pn 20}%
\special{pa 2237 480}%
\special{pa 3167 480}%
\special{fp}%
%
\special{pn 20}%
\special{pa 2621 1118}%
\special{pa 3327 1434}%
\special{fp}%
%
\special{pn 20}%
\special{pa 1684 1123}%
\special{pa 2391 1439}%
\special{fp}%
%
\special{pn 20}%
\special{pa 2395 1441}%
\special{pa 3324 1441}%
\special{fp}%
%
\special{pn 20}%
\special{pa 1933 1230}%
\special{pa 2862 1230}%
\special{fp}%
%
\special{pn 20}%
\special{pa 2195 1349}%
\special{pa 3124 1349}%
\special{fp}%
%
\special{pn 20}%
\special{pa 3337 400}%
\special{pa 3831 966}%
\special{fp}%
\special{pa 3324 1439}%
\special{pa 3828 966}%
\special{fp}%
%
\special{pn 20}%
\special{pa 2122 945}%
\special{pa 2187 945}%
\special{fp}%
\special{pa 2187 945}%
\special{pa 2119 986}%
\special{fp}%
\special{pa 2119 986}%
\special{pa 2156 908}%
\special{fp}%
\special{pa 2156 908}%
\special{pa 2175 989}%
\special{fp}%
\special{pa 2175 989}%
\special{pa 2126 947}%
\special{fp}%
%
\special{pn 20}%
\special{sh 0.600}%
\special{ar 1684 952 45 34  0.0000000 6.2831853}%
%
\special{pn 20}%
\special{pa 2130 963}%
\special{pa 1974 957}%
\special{fp}%
\special{sh 1}%
\special{pa 1974 957}%
\special{pa 2040 980}%
\special{pa 2027 959}%
\special{pa 2041 940}%
\special{pa 1974 957}%
\special{fp}%
\special{pa 2149 957}%
\special{pa 2149 1077}%
\special{fp}%
\special{sh 1}%
\special{pa 2149 1077}%
\special{pa 2169 1010}%
\special{pa 2149 1024}%
\special{pa 2129 1010}%
\special{pa 2149 1077}%
\special{fp}%
\put(14.1000,-9.7000){\makebox(0,0)[lb]{$P_4$}}%
\put(23.1700,-9.5700){\makebox(0,0)[lb]{$P_7$}}%
\put(19.1800,-9.4700){\makebox(0,0)[lb]{$l_1$}}%
\put(21.6300,-10.7700){\makebox(0,0)[lb]{$m_1$}}%
\put(14.3000,-6.8900){\makebox(0,0)[lb]{$C_6$}}%
%
\special{pn 20}%
\special{pa 1728 956}%
\special{pa 1884 956}%
\special{fp}%
\special{sh 1}%
\special{pa 1884 956}%
\special{pa 1817 936}%
\special{pa 1831 956}%
\special{pa 1817 976}%
\special{pa 1884 956}%
\special{fp}%
%
\special{pn 20}%
\special{pa 1681 956}%
\special{pa 1688 1120}%
\special{fp}%
\special{sh 1}%
\special{pa 1688 1120}%
\special{pa 1705 1053}%
\special{pa 1686 1067}%
\special{pa 1665 1054}%
\special{pa 1688 1120}%
\special{fp}%
\put(17.4200,-9.4700){\makebox(0,0)[lb]{$L_1$}}%
\put(14.5000,-11.0500){\makebox(0,0)[lb]{$M_1$}}%
%
\special{pn 8}%
\special{pa 1572 737}%
\special{pa 1681 807}%
\special{fp}%
\special{sh 1}%
\special{pa 1681 807}%
\special{pa 1636 754}%
\special{pa 1636 778}%
\special{pa 1614 788}%
\special{pa 1681 807}%
\special{fp}%
%
\special{pn 20}%
\special{pa 4852 833}%
\special{pa 4852 1190}%
\special{fp}%
%
\special{pn 20}%
\special{pa 4852 840}%
\special{pa 6076 840}%
\special{fp}%
%
\special{pn 20}%
\special{pa 6076 849}%
\special{pa 6076 1206}%
\special{fp}%
%
\special{pn 20}%
\special{pa 4852 1198}%
\special{pa 6076 1198}%
\special{fp}%
%
\special{pn 20}%
\special{sh 0.600}%
\special{ar 4852 1039 45 35  0.0000000 6.2831853}%
%
\special{pn 20}%
\special{pa 5514 1055}%
\special{pa 5580 1055}%
\special{fp}%
\special{pa 5580 1055}%
\special{pa 5511 1097}%
\special{fp}%
\special{pa 5511 1097}%
\special{pa 5548 1020}%
\special{fp}%
\special{pa 5548 1020}%
\special{pa 5568 1100}%
\special{fp}%
\special{pa 5568 1100}%
\special{pa 5519 1057}%
\special{fp}%
%
\special{pn 20}%
\special{pa 4852 1047}%
\special{pa 4858 1212}%
\special{fp}%
\special{sh 1}%
\special{pa 4858 1212}%
\special{pa 4876 1145}%
\special{pa 4856 1159}%
\special{pa 4836 1146}%
\special{pa 4858 1212}%
\special{fp}%
%
\special{pn 20}%
\special{pa 4868 1047}%
\special{pa 5024 1047}%
\special{fp}%
\special{sh 1}%
\special{pa 5024 1047}%
\special{pa 4957 1027}%
\special{pa 4971 1047}%
\special{pa 4957 1067}%
\special{pa 5024 1047}%
\special{fp}%
%
\special{pn 20}%
\special{pa 5523 1063}%
\special{pa 5367 1057}%
\special{fp}%
\special{sh 1}%
\special{pa 5367 1057}%
\special{pa 5433 1080}%
\special{pa 5420 1059}%
\special{pa 5434 1040}%
\special{pa 5367 1057}%
\special{fp}%
\special{pa 5543 1057}%
\special{pa 5543 1177}%
\special{fp}%
\special{sh 1}%
\special{pa 5543 1177}%
\special{pa 5563 1110}%
\special{pa 5543 1124}%
\special{pa 5523 1110}%
\special{pa 5543 1177}%
\special{fp}%
%
\special{pn 8}%
\special{pa 4716 833}%
\special{pa 4824 902}%
\special{fp}%
\special{sh 1}%
\special{pa 4824 902}%
\special{pa 4779 849}%
\special{pa 4779 873}%
\special{pa 4757 883}%
\special{pa 4824 902}%
\special{fp}%
\put(44.6000,-8.1600){\makebox(0,0)[lb]{$C_6$}}%
%
\special{pn 8}%
\special{pa 4852 1047}%
\special{pa 6067 1047}%
\special{dt 0.045}%
\special{pa 6067 1047}%
\special{pa 6066 1047}%
\special{dt 0.045}%
%
\special{pn 8}%
\special{pa 4868 968}%
\special{pa 6058 968}%
\special{dt 0.045}%
\special{pa 6058 968}%
\special{pa 6057 968}%
\special{dt 0.045}%
%
\special{pn 8}%
\special{pa 4868 912}%
\special{pa 6076 912}%
\special{dt 0.045}%
\special{pa 6076 912}%
\special{pa 6075 912}%
\special{dt 0.045}%
%
\special{pn 8}%
\special{pa 4852 1119}%
\special{pa 6067 1119}%
\special{dt 0.045}%
\special{pa 6067 1119}%
\special{pa 6066 1119}%
\special{dt 0.045}%
\put(39.0000,-10.9400){\makebox(0,0)[lb]{$\supset$}}%
%
\special{pn 20}%
\special{pa 1759 1897}%
\special{pa 1759 2404}%
\special{fp}%
\special{pa 1759 1897}%
\special{pa 2442 1780}%
\special{fp}%
\special{pa 2442 1774}%
\special{pa 2966 2164}%
\special{fp}%
%
\special{pn 20}%
\special{pa 2442 1774}%
\special{pa 2442 2248}%
\special{fp}%
\special{pa 1759 2404}%
\special{pa 2442 2248}%
\special{fp}%
\special{pa 2442 2248}%
\special{pa 2966 2632}%
\special{fp}%
%
\special{pn 20}%
\special{pa 2966 2164}%
\special{pa 2966 2632}%
\special{fp}%
%
\special{pn 20}%
\special{pa 3561 1949}%
\special{pa 3561 2488}%
\special{fp}%
%
\special{pn 20}%
\special{pa 3561 1940}%
\special{pa 4708 1940}%
\special{fp}%
%
\special{pn 20}%
\special{pa 3570 2481}%
\special{pa 4717 2481}%
\special{fp}%
%
\special{pn 20}%
\special{pa 4717 1940}%
\special{pa 4717 2481}%
\special{fp}%
%
\special{pn 20}%
\special{pa 3119 2234}%
\special{pa 3425 2234}%
\special{fp}%
\special{sh 1}%
\special{pa 3425 2234}%
\special{pa 3358 2214}%
\special{pa 3372 2234}%
\special{pa 3358 2254}%
\special{pa 3425 2234}%
\special{fp}%
%
\special{pn 20}%
\special{pa 1623 2210}%
\special{pa 1130 2210}%
\special{fp}%
\special{sh 1}%
\special{pa 1130 2210}%
\special{pa 1197 2230}%
\special{pa 1183 2210}%
\special{pa 1197 2190}%
\special{pa 1130 2210}%
\special{fp}%
%
\special{pn 20}%
\special{pa 2660 2171}%
\special{pa 2725 2171}%
\special{fp}%
\special{pa 2725 2171}%
\special{pa 2657 2213}%
\special{fp}%
\special{pa 2657 2213}%
\special{pa 2694 2136}%
\special{fp}%
\special{pa 2694 2136}%
\special{pa 2714 2216}%
\special{fp}%
\special{pa 2714 2216}%
\special{pa 2664 2173}%
\special{fp}%
%
\special{pn 20}%
\special{pa 4692 2234}%
\special{pa 4757 2234}%
\special{fp}%
\special{pa 4757 2234}%
\special{pa 4688 2276}%
\special{fp}%
\special{pa 4688 2276}%
\special{pa 4726 2199}%
\special{fp}%
\special{pa 4726 2199}%
\special{pa 4745 2279}%
\special{fp}%
\special{pa 4745 2279}%
\special{pa 4696 2237}%
\special{fp}%
%
\special{pn 20}%
\special{pa 5465 1920}%
\special{pa 5465 2428}%
\special{fp}%
\special{pa 5465 1920}%
\special{pa 6148 1804}%
\special{fp}%
\special{pa 6148 1797}%
\special{pa 6672 2188}%
\special{fp}%
%
\special{pn 20}%
\special{pa 6148 1797}%
\special{pa 6148 2272}%
\special{fp}%
\special{pa 5465 2428}%
\special{pa 6148 2272}%
\special{fp}%
\special{pa 6148 2272}%
\special{pa 6672 2656}%
\special{fp}%
%
\special{pn 20}%
\special{pa 6672 2188}%
\special{pa 6672 2656}%
\special{fp}%
%
\special{pn 20}%
\special{pa 5329 2234}%
\special{pa 4836 2234}%
\special{fp}%
\special{sh 1}%
\special{pa 4836 2234}%
\special{pa 4903 2254}%
\special{pa 4889 2234}%
\special{pa 4903 2214}%
\special{pa 4836 2234}%
\special{fp}%
%
\special{pn 20}%
\special{pa 6646 2401}%
\special{pa 6712 2401}%
\special{fp}%
\special{pa 6712 2401}%
\special{pa 6643 2443}%
\special{fp}%
\special{pa 6643 2443}%
\special{pa 6680 2366}%
\special{fp}%
\special{pa 6680 2366}%
\special{pa 6700 2446}%
\special{fp}%
\special{pa 6700 2446}%
\special{pa 6651 2403}%
\special{fp}%
%
\special{pn 20}%
\special{pa 1589 2862}%
\special{pa 1589 3403}%
\special{fp}%
%
\special{pn 20}%
\special{pa 1589 2854}%
\special{pa 2736 2854}%
\special{fp}%
%
\special{pn 20}%
\special{pa 1598 3395}%
\special{pa 2745 3395}%
\special{fp}%
%
\special{pn 20}%
\special{pa 2745 2854}%
\special{pa 2745 3395}%
\special{fp}%
%
\special{pn 20}%
\special{pa 1147 3148}%
\special{pa 1453 3148}%
\special{fp}%
\special{sh 1}%
\special{pa 1453 3148}%
\special{pa 1386 3128}%
\special{pa 1400 3148}%
\special{pa 1386 3168}%
\special{pa 1453 3148}%
\special{fp}%
%
\special{pn 20}%
\special{pa 2720 3148}%
\special{pa 2785 3148}%
\special{fp}%
\special{pa 2785 3148}%
\special{pa 2716 3191}%
\special{fp}%
\special{pa 2716 3191}%
\special{pa 2754 3113}%
\special{fp}%
\special{pa 2754 3113}%
\special{pa 2773 3193}%
\special{fp}%
\special{pa 2773 3193}%
\special{pa 2724 3151}%
\special{fp}%
%
\special{pn 20}%
\special{sh 0.600}%
\special{ar 1750 2171 46 34  0.0000000 6.2831853}%
%
\special{pn 20}%
\special{sh 0.600}%
\special{ar 3552 2250 46 35  0.0000000 6.2831853}%
%
\special{pn 20}%
\special{sh 0.600}%
\special{ar 5465 2243 45 34  0.0000000 6.2831853}%
%
\special{pn 20}%
\special{sh 0.600}%
\special{ar 1598 3157 45 34  0.0000000 6.2831853}%
%
\special{pn 20}%
\special{pa 3510 2851}%
\special{pa 3510 3358}%
\special{fp}%
\special{pa 3510 2851}%
\special{pa 4193 2734}%
\special{fp}%
\special{pa 4193 2727}%
\special{pa 4717 3117}%
\special{fp}%
%
\special{pn 20}%
\special{pa 4193 2727}%
\special{pa 4193 3202}%
\special{fp}%
\special{pa 3510 3358}%
\special{pa 4193 3202}%
\special{fp}%
\special{pa 4193 3202}%
\special{pa 4717 3586}%
\special{fp}%
%
\special{pn 20}%
\special{pa 4717 3117}%
\special{pa 4717 3586}%
\special{fp}%
%
\special{pn 20}%
\special{pa 3374 3164}%
\special{pa 2881 3164}%
\special{fp}%
\special{sh 1}%
\special{pa 2881 3164}%
\special{pa 2948 3184}%
\special{pa 2934 3164}%
\special{pa 2948 3144}%
\special{pa 2881 3164}%
\special{fp}%
%
\special{pn 20}%
\special{pa 4692 3331}%
\special{pa 4757 3331}%
\special{fp}%
\special{pa 4757 3331}%
\special{pa 4688 3373}%
\special{fp}%
\special{pa 4688 3373}%
\special{pa 4726 3296}%
\special{fp}%
\special{pa 4726 3296}%
\special{pa 4745 3376}%
\special{fp}%
\special{pa 4745 3376}%
\special{pa 4696 3333}%
\special{fp}%
%
\special{pn 20}%
\special{sh 0.600}%
\special{ar 3510 3172 45 35  0.0000000 6.2831853}%
%
\special{pn 20}%
\special{pa 5286 2870}%
\special{pa 5286 3411}%
\special{fp}%
%
\special{pn 20}%
\special{pa 5286 2862}%
\special{pa 6434 2862}%
\special{fp}%
%
\special{pn 20}%
\special{pa 5295 3403}%
\special{pa 6442 3403}%
\special{fp}%
%
\special{pn 20}%
\special{pa 6442 2862}%
\special{pa 6442 3403}%
\special{fp}%
%
\special{pn 20}%
\special{pa 4844 3157}%
\special{pa 5150 3157}%
\special{fp}%
\special{sh 1}%
\special{pa 5150 3157}%
\special{pa 5083 3137}%
\special{pa 5097 3157}%
\special{pa 5083 3177}%
\special{pa 5150 3157}%
\special{fp}%
%
\special{pn 20}%
\special{pa 6417 3157}%
\special{pa 6482 3157}%
\special{fp}%
\special{pa 6482 3157}%
\special{pa 6414 3198}%
\special{fp}%
\special{pa 6414 3198}%
\special{pa 6451 3121}%
\special{fp}%
\special{pa 6451 3121}%
\special{pa 6471 3201}%
\special{fp}%
\special{pa 6471 3201}%
\special{pa 6421 3158}%
\special{fp}%
%
\special{pn 20}%
\special{sh 0.600}%
\special{ar 5295 3164 45 35  0.0000000 6.2831853}%
%
\special{pn 8}%
\special{pa 5304 3164}%
\special{pa 6434 3164}%
\special{dt 0.045}%
\special{pa 6434 3164}%
\special{pa 6433 3164}%
\special{dt 0.045}%
%
\special{pn 20}%
\special{pa 1415 3962}%
\special{pa 959 3962}%
\special{fp}%
\special{sh 1}%
\special{pa 959 3962}%
\special{pa 1026 3982}%
\special{pa 1012 3962}%
\special{pa 1026 3942}%
\special{pa 959 3962}%
\special{fp}%
%
\special{pn 20}%
\special{pa 3207 3962}%
\special{pa 2751 3962}%
\special{fp}%
\special{sh 1}%
\special{pa 2751 3962}%
\special{pa 2818 3982}%
\special{pa 2804 3962}%
\special{pa 2818 3942}%
\special{pa 2751 3962}%
\special{fp}%
%
\special{pn 20}%
\special{pa 5040 3977}%
\special{pa 4583 3977}%
\special{fp}%
\special{sh 1}%
\special{pa 4583 3977}%
\special{pa 4650 3997}%
\special{pa 4636 3977}%
\special{pa 4650 3957}%
\special{pa 4583 3977}%
\special{fp}%
%
\special{pn 20}%
\special{pa 1488 3681}%
\special{pa 1488 4199}%
\special{fp}%
%
\special{pn 20}%
\special{pa 1488 3672}%
\special{pa 2550 3672}%
\special{fp}%
%
\special{pn 20}%
\special{pa 1496 4190}%
\special{pa 2557 4190}%
\special{fp}%
%
\special{pn 20}%
\special{pa 2557 3672}%
\special{pa 2557 4190}%
\special{fp}%
%
\special{pn 20}%
\special{pa 2657 4122}%
\special{pa 2717 4122}%
\special{fp}%
\special{pa 2717 4122}%
\special{pa 2654 4161}%
\special{fp}%
\special{pa 2654 4161}%
\special{pa 2688 4088}%
\special{fp}%
\special{pa 2688 4088}%
\special{pa 2707 4165}%
\special{fp}%
\special{pa 2707 4165}%
\special{pa 2661 4124}%
\special{fp}%
%
\special{pn 20}%
\special{sh 0.600}%
\special{ar 1737 4114 42 34  0.0000000 6.2831853}%
%
\special{pn 8}%
\special{pa 1503 3962}%
\special{pa 2550 3962}%
\special{dt 0.045}%
\special{pa 2550 3962}%
\special{pa 2549 3962}%
\special{dt 0.045}%
%
\special{pn 8}%
\special{pa 1486 3977}%
\special{pa 1521 3975}%
\special{pa 1555 3975}%
\special{pa 1587 3978}%
\special{pa 1617 3985}%
\special{pa 1643 3998}%
\special{pa 1665 4017}%
\special{pa 1683 4042}%
\special{pa 1700 4071}%
\special{pa 1714 4102}%
\special{pa 1728 4135}%
\special{pa 1729 4137}%
\special{sp -0.045}%
%
\special{pn 8}%
\special{pa 1729 4122}%
\special{pa 2538 4122}%
\special{dt 0.045}%
\special{pa 2538 4122}%
\special{pa 2537 4122}%
\special{dt 0.045}%
%
\special{pn 8}%
\special{pa 2555 4122}%
\special{pa 2681 4122}%
\special{dt 0.045}%
\special{pa 2681 4122}%
\special{pa 2680 4122}%
\special{dt 0.045}%
%
\special{pn 20}%
\special{pa 2547 3970}%
\special{pa 2583 3977}%
\special{pa 2616 3985}%
\special{pa 2644 3998}%
\special{pa 2665 4017}%
\special{pa 2678 4043}%
\special{pa 2685 4075}%
\special{pa 2687 4110}%
\special{pa 2688 4122}%
\special{sp}%
%
\special{pn 20}%
\special{pa 3359 3696}%
\special{pa 3359 4213}%
\special{fp}%
%
\special{pn 20}%
\special{pa 3359 3688}%
\special{pa 4421 3688}%
\special{fp}%
%
\special{pn 20}%
\special{pa 3367 4205}%
\special{pa 4428 4205}%
\special{fp}%
%
\special{pn 20}%
\special{pa 4428 3688}%
\special{pa 4428 4205}%
\special{fp}%
%
\special{pn 20}%
\special{pa 4740 4296}%
\special{pa 4801 4296}%
\special{fp}%
\special{pa 4801 4296}%
\special{pa 4737 4337}%
\special{fp}%
\special{pa 4737 4337}%
\special{pa 4771 4263}%
\special{fp}%
\special{pa 4771 4263}%
\special{pa 4790 4340}%
\special{fp}%
\special{pa 4790 4340}%
\special{pa 4745 4299}%
\special{fp}%
%
\special{pn 20}%
\special{sh 0.600}%
\special{ar 3671 4313 42 32  0.0000000 6.2831853}%
%
\special{pn 8}%
\special{pa 3375 3977}%
\special{pa 4421 3977}%
\special{dt 0.045}%
\special{pa 4421 3977}%
\special{pa 4420 3977}%
\special{dt 0.045}%
%
\special{pn 8}%
\special{pa 3356 3985}%
\special{pa 3391 3982}%
\special{pa 3424 3984}%
\special{pa 3454 3995}%
\special{pa 3479 4019}%
\special{pa 3499 4061}%
\special{pa 3512 4111}%
\special{pa 3519 4146}%
\special{pa 3522 4145}%
\special{pa 3522 4137}%
\special{sp -0.045}%
%
\special{pn 8}%
\special{pa 3513 4122}%
\special{pa 3543 4136}%
\special{pa 3572 4152}%
\special{pa 3597 4171}%
\special{pa 3617 4194}%
\special{pa 3632 4221}%
\special{pa 3644 4251}%
\special{pa 3653 4284}%
\special{pa 3656 4296}%
\special{sp -0.045}%
%
\special{pn 8}%
\special{pa 3513 4122}%
\special{pa 4505 4122}%
\special{dt 0.045}%
\special{pa 4505 4122}%
\special{pa 4504 4122}%
\special{dt 0.045}%
%
\special{pn 20}%
\special{pa 4426 3977}%
\special{pa 4460 3985}%
\special{pa 4492 3994}%
\special{pa 4521 4006}%
\special{pa 4545 4024}%
\special{pa 4564 4047}%
\special{pa 4578 4075}%
\special{pa 4589 4107}%
\special{pa 4598 4140}%
\special{pa 4599 4145}%
\special{sp}%
%
\special{pn 20}%
\special{pa 4591 4122}%
\special{pa 4624 4120}%
\special{pa 4656 4123}%
\special{pa 4686 4135}%
\special{pa 4713 4160}%
\special{pa 4735 4203}%
\special{pa 4752 4252}%
\special{pa 4764 4291}%
\special{pa 4770 4305}%
\special{pa 4771 4290}%
\special{sp}%
%
\special{pn 20}%
\special{pa 3671 4313}%
\special{pa 4780 4313}%
\special{fp}%
%
\special{pn 20}%
\special{pa 4474 4130}%
\special{pa 4599 4130}%
\special{fp}%
%
\special{pn 8}%
\special{pa 5113 3688}%
\special{pa 5113 4205}%
\special{dt 0.045}%
\special{pa 5113 4205}%
\special{pa 5113 4204}%
\special{dt 0.045}%
%
\special{pn 8}%
\special{pa 5113 3681}%
\special{pa 6174 3681}%
\special{dt 0.045}%
\special{pa 6174 3681}%
\special{pa 6173 3681}%
\special{dt 0.045}%
%
\special{pn 8}%
\special{pa 5120 4199}%
\special{pa 6182 4199}%
\special{dt 0.045}%
\special{pa 6182 4199}%
\special{pa 6181 4199}%
\special{dt 0.045}%
%
\special{pn 8}%
\special{pa 6182 3681}%
\special{pa 6182 4199}%
\special{dt 0.045}%
\special{pa 6182 4199}%
\special{pa 6182 4198}%
\special{dt 0.045}%
%
\special{pn 8}%
\special{pa 5129 3970}%
\special{pa 6174 3970}%
\special{dt 0.045}%
\special{pa 6174 3970}%
\special{pa 6173 3970}%
\special{dt 0.045}%
%
\special{pn 8}%
\special{pa 5110 3977}%
\special{pa 5145 3974}%
\special{pa 5178 3976}%
\special{pa 5208 3988}%
\special{pa 5233 4012}%
\special{pa 5252 4054}%
\special{pa 5265 4104}%
\special{pa 5273 4139}%
\special{pa 5275 4136}%
\special{pa 5275 4130}%
\special{sp -0.045}%
%
\special{pn 8}%
\special{pa 5267 4114}%
\special{pa 5297 4129}%
\special{pa 5325 4145}%
\special{pa 5350 4164}%
\special{pa 5370 4187}%
\special{pa 5385 4214}%
\special{pa 5397 4245}%
\special{pa 5406 4277}%
\special{pa 5409 4290}%
\special{sp -0.045}%
%
\special{pn 8}%
\special{pa 5267 4114}%
\special{pa 6258 4114}%
\special{dt 0.045}%
\special{pa 6258 4114}%
\special{pa 6257 4114}%
\special{dt 0.045}%
%
\special{pn 8}%
\special{pa 6179 3970}%
\special{pa 6213 3977}%
\special{pa 6245 3986}%
\special{pa 6274 3998}%
\special{pa 6299 4015}%
\special{pa 6317 4039}%
\special{pa 6331 4067}%
\special{pa 6342 4099}%
\special{pa 6351 4132}%
\special{pa 6352 4137}%
\special{sp -0.045}%
%
\special{pn 8}%
\special{pa 6344 4114}%
\special{pa 6377 4113}%
\special{pa 6409 4116}%
\special{pa 6439 4128}%
\special{pa 6465 4153}%
\special{pa 6488 4195}%
\special{pa 6505 4243}%
\special{pa 6517 4283}%
\special{pa 6524 4298}%
\special{pa 6525 4281}%
\special{sp -0.045}%
%
\special{pn 8}%
\special{pa 5425 4305}%
\special{pa 6533 4305}%
\special{dt 0.045}%
\special{pa 6533 4305}%
\special{pa 6532 4305}%
\special{dt 0.045}%
%
\special{pn 8}%
\special{pa 6227 4122}%
\special{pa 6352 4122}%
\special{dt 0.045}%
\special{pa 6352 4122}%
\special{pa 6351 4122}%
\special{dt 0.045}%
%
\special{pn 20}%
\special{pa 5425 4296}%
\special{pa 5455 4311}%
\special{pa 5484 4327}%
\special{pa 5511 4344}%
\special{pa 5535 4364}%
\special{pa 5555 4388}%
\special{pa 5570 4415}%
\special{pa 5582 4445}%
\special{pa 5593 4476}%
\special{pa 5598 4495}%
\special{sp}%
%
\special{pn 20}%
\special{pa 6518 4313}%
\special{pa 6549 4327}%
\special{pa 6578 4343}%
\special{pa 6605 4360}%
\special{pa 6629 4380}%
\special{pa 6649 4403}%
\special{pa 6664 4430}%
\special{pa 6676 4460}%
\special{pa 6685 4492}%
\special{pa 6690 4510}%
\special{sp}%
%
\special{pn 20}%
\special{pa 5589 4503}%
\special{pa 6675 4503}%
\special{fp}%
\put(10.9000,-21.6000){\makebox(0,0)[lb]{Step 6}}%
\put(30.0000,-22.1000){\makebox(0,0)[lb]{Step 6}}%
\put(48.4000,-22.0200){\makebox(0,0)[lb]{Step 7}}%
\put(10.0000,-31.1200){\makebox(0,0)[lb]{Step 7}}%
\put(28.8000,-31.2800){\makebox(0,0)[lb]{Step 8}}%
\put(47.4000,-31.2800){\makebox(0,0)[lb]{Step 8}}%
\put(9.4000,-39.1700){\makebox(0,0)[lb]{Step 9}}%
\put(27.3500,-39.2500){\makebox(0,0)[lb]{Step 10}}%
\put(45.0100,-39.5000){\makebox(0,0)[lb]{Step 11}}%
%
\special{pn 8}%
\special{pa 3456 1955}%
\special{pa 3564 2024}%
\special{fp}%
\special{sh 1}%
\special{pa 3564 2024}%
\special{pa 3519 1971}%
\special{pa 3519 1995}%
\special{pa 3497 2005}%
\special{pa 3564 2024}%
\special{fp}%
\put(32.0000,-19.3800){\makebox(0,0)[lb]{$C_7$}}%
%
\special{pn 8}%
\special{pa 1496 2881}%
\special{pa 1604 2951}%
\special{fp}%
\special{sh 1}%
\special{pa 1604 2951}%
\special{pa 1559 2898}%
\special{pa 1559 2922}%
\special{pa 1537 2932}%
\special{pa 1604 2951}%
\special{fp}%
\put(12.4000,-28.6500){\makebox(0,0)[lb]{$C_8$}}%
\put(56.1000,-31.7100){\makebox(0,0)[lb]{$C_9$}}%
\put(19.7600,-41.4500){\makebox(0,0)[lb]{$C_{10}$}}%
\put(38.7200,-44.4600){\makebox(0,0)[lb]{$C_{11}$}}%
\put(25.8000,-23.6400){\makebox(0,0)[lb]{$S_1$}}%
\put(62.0000,-23.1200){\makebox(0,0)[lb]{$S_2$}}%
\put(42.2000,-32.3000){\makebox(0,0)[lb]{$S_3$}}%
%
\special{pn 20}%
\special{pa 1760 2170}%
\special{pa 2090 2090}%
\special{fp}%
\special{sh 1}%
\special{pa 2090 2090}%
\special{pa 2020 2086}%
\special{pa 2038 2103}%
\special{pa 2030 2125}%
\special{pa 2090 2090}%
\special{fp}%
%
\special{pn 20}%
\special{pa 5470 2230}%
\special{pa 5810 2150}%
\special{fp}%
\special{sh 1}%
\special{pa 5810 2150}%
\special{pa 5741 2146}%
\special{pa 5758 2162}%
\special{pa 5750 2185}%
\special{pa 5810 2150}%
\special{fp}%
%
\special{pn 20}%
\special{pa 3510 3170}%
\special{pa 3880 3080}%
\special{fp}%
\special{sh 1}%
\special{pa 3880 3080}%
\special{pa 3810 3076}%
\special{pa 3828 3093}%
\special{pa 3820 3115}%
\special{pa 3880 3080}%
\special{fp}%
\put(18.5000,-21.0000){\makebox(0,0)[lb]{$L_2$}}%
\put(55.4000,-21.5000){\makebox(0,0)[lb]{$L_3$}}%
\put(35.9000,-30.8000){\makebox(0,0)[lb]{$L_4$}}%
\end{picture}%
\label{fig:third2}
\caption{}
\end{figure}

\vspace{0.3cm}

{\bf 4.3. Resolution of accessible singular points $P_4$ and $P_7$ }

Let us start by summarizing the steps which are needed to resolve the accessible singular points of $\tilde v$ at $P_4$ and $P_7$:

$$ P_4=\{(u,v,w)=(0,-1,0)\}, \  P_7=\{(u,v,w)=(-3,-1,0)\}.$$

\begin{enumerate}
\item{\bf Step 6}: We will blow up along the curve $C_6$ and blow down the surface $S_1$ (see Figure 4).
\item{\bf Step 7}: We will blow up along the curve $C_7$ and blow down the surface $S_2$ (see Figure 4).
\item{\bf Step 8}: We will blow up along the curve $C_{8}$ and blow down the surface $S_3$  (see Figure 4).
\item{\bf Step 9}: We will blow up along the curve $C_{9}$ (see Figure 4).
\item{\bf Step 10}: We will blow up along the curve $C_{10}$ (see Figure 4).
\item{\bf Step 11}: We will blow up along the curve $C_{11}$ (see Figure 4).
\end{enumerate}

In a neighborhood of $P_7$, the system is written as

\begin{equation*}
  \left\{
  \begin{aligned}
   \frac {dl_1}{dt} &=\frac{3l_1-(l_1)^2-6m_1+2l_1m_1}{n_1}+n_1\alpha_1
   +3n_1\alpha_3
   -l_1n_1\alpha_3\\
   \frac {dm_1}{dt} &=\frac{-3m_1+3(m_1)^2}{n_1}+2t-2tm_1+n_1\alpha_2
   +n_1\alpha_3-m_1n_1\alpha_3\\
   \frac {dn_1}{dt} &=-1+2m_1-tn_1-\alpha_3(n_1)^2.
   \end{aligned}
  \right. 
\end{equation*}

In a neighborhood of $P_4$, the system is written as
\begin{equation*}
  \left\{
  \begin{aligned}
   \frac {dL_1}{dt} &=\frac{-3L_1-(L_1)^2+2L_1M_1}{N_1}+\alpha_1N_1
   -\alpha_3L_1N_1\\
   \frac {dM_1}{dt} &=2t-2tM_1+\frac{-3M_1+3(M_1)^2}{N_1}+\alpha_2N_1
   +\alpha_3N_1-\alpha_3M_1N_1\\
   \frac {dN_1}{dt} &=-1+2M_1-tN_1-\alpha_3(N_1)^2.
   \end{aligned}
  \right. 
\end{equation*}

{\bf 4.3.1.} We blow up along the curve $C_6=\{(L_1, M_1, N_1)| L_1=N_1=0\} \subset {\Bbb P}^1$ and blow down the surface $S_1=\{(l_1,m_1,n_1)|n_1=0\}$

$$
L_2=\frac{L_1}{N_1} \;, \;\;\; M_2=M_1 \;, \;\;\; N_2=N_1,
$$
$$
l_2=\frac{n_1}{l_1-3} \;, \;\;\;  m_2=m_1 \;, \;\;\; n_2=n_1.
$$

In a neighborhood of  $\{(l_2,m_2.n_2)=(0,0,0)\}$, the system is written as

\begin{equation*}
  \left\{
  \begin{aligned}
   \frac {dl_2}{dt} &=1-tl_2+\frac{2l_2}{n_2}-\alpha_1(l_2)^2\\
   \frac {dm_2}{dt} &=2t-2tm_2-\frac{3m_2-3(m_2)^2}{n_2}+\alpha_2n_2
   +\alpha_3n_2-\alpha_3m_2n_2\\
   \frac {dn_2}{dt} &=-1+2m_2-tn_2-\alpha_3(n_2)^2.
   \end{aligned}
  \right. 
\end{equation*}

{\bf 4.3.2.} We blow up along the curve $C_7=\{(L_2, M_2, N_2)| L_2=N_2=0\} \subset {\Bbb P}^1$ and blow down the surface $S_2=\{(l_2,m_2,n_2)|n_2=0\}$

$$
L_3=\frac{L_2}{N_2} \;, \;\;\; M_3=M_2 \;, \;\;\; N_3=N_2,
$$
$$
l_3=l_2n_2 \;, \;\;\; m_3=m_2 \;, \;\;\; n_3=n_2.
$$

In a neighborhood of  $\{(l_3,m_3,n_3)=(0,0,0)\}$, the system is written as

\begin{equation*}
  \left\{
  \begin{aligned}
   \frac {dl_3}{dt} &=-2tl_3+\frac{l_3+2l_3m_3-\alpha_2(l_3)^2}{n_3}
   +n_3-\alpha_3l_3n_3\\
   \frac {dm_3}{dt} &=2t-2tm_3-\frac{3m_3-3(m_3)^2}{n_3}+\alpha_2n_3
   +\alpha_3n_3-\alpha_3m_3n_3\\
   \frac {dn_3}{dt} &=-1+2m_3-tn_3-\alpha_3(n_3)^2.
   \end{aligned}
  \right. 
\end{equation*}

{\bf 4.3.3.} We blow up along the curve $C_{8}=\{(L_3, M_3, N_3)| L_3=N_3=0\} \subset {\Bbb P}^1$ and blow down the surface $S_3=\{(l_3,m_3,n_3)|n_3=0\}$

$$
L_4=\frac{L_3-\alpha_1}{N_3} \;, \;\;\; M_4=M_3 \;, \;\;\; N_4=N_3,
$$
$$
l_4=\frac{l_3n_3}{1-\alpha_1l_3} \;, \;\;\; m_4=m_3 \;, \;\;\; n_4=n_3.
$$

{\bf 4.3.4.} We blow up along the curve $C_{9}=\{(L_4, M_4, N_4)| M_4=N_4=0\} \cup \{(l_4, m_4, n_4)| m_4=n_4=0\} \cong {\Bbb P}^1$

$$
L_5=L_4 \;, \;\;\; M_5=\frac{M_4}{N_4} \;, \;\;\; N_5=N_4,
$$
$$
l_5=l_4 \;, \;\;\; m_5=\frac{m_4}{n_4} \;, \;\;\; n_5=n_4.
$$

{\bf 4.3.5.} We blow up along the curve $C_{10}=\{(L_5, M_5, N_5)| M_5-t=N_5=0\} \cup \{(l_5, m_5, n_5)| m_5-t=n_5=0\} \cong {\Bbb P}^1$

$$
L_6=L_5 \;, \;\;\; M_6=\frac{M_5-t}{N_5} \;, \;\;\; N_6=N_5,
$$
$$
l_6=l_5 \;, \;\;\; m_6=\frac{m_5-t}{n_5} \;, \;\;\; n_6=n_5.
$$

{\bf 4.3.6.} We blow up along the curve $C_{11}=\{(L_6, M_6, N_6)| M_6+1-\alpha_2-\alpha_3=N_6=0\} \cup \{(l_6, m_6, n_6)| m_6+1-\alpha_2-\alpha_3=n_6=0\} \cong {\Bbb P}^1$

$$
x_4=L_6 \;, \;\;\; y_4=\frac{M_6+1-\alpha_2-\alpha_3}{N_6} \;, \;\;\; z_4=N_6,
$$
$$
x_7=l_6 \;, \;\;\; y_7=\frac{m_6+1-\alpha_2-\alpha_3}{n_6} \;, \;\;\; z_7=n_6.
$$
We have resolved the accessible singular points $P_4$ and $P_7$.

\vspace{0.3cm}

{\it Acknowledgement.} The author thanks Professors M.-H. Saito, K. Takano and Y. Yamada for giving helpful suggestions and encouragement, and Professor W. Rossman for checking English.


\begin{thebibliography}{99}




\bibitem[1]{1} K. Okamoto, {\em Sur les 
feuilletages associ\'es aux \'equations du second ordre \`a points critiques fixes de P. Painlev\'e, Espaces des conditions initiales}, Japan. 
J. Math., {\bf 5}, 1979, 1--79.  


\bibitem[2]{2} M. Noumi and Y. Yamada, 
        {\em Higher order Painlev\`e equations of type ${A_l}^{(1)}$}, Funkcial. Ekvac., {\bf 41}, 1998, 483--503.


\bibitem[3]{3} T. Shioda and K. Takano, On some Hamiltonian structures of Painlev\'e systems I, Funkcial. Ekvac., {\bf 40}, 1997, 271--291. 



\bibitem[4]{4} Y. Sasano, Coupled Painlev\'e II systems in dimension 4 and the systems of type ${A_4}^{(1)}$, preprint.



\bibitem[5]{5} N. Tahara, An augmentation of the phase space of the system of type ${A_4}^{(1)}$, Kyushu J. Math. {\bf 58}, 2004, 393--425. 


\end{thebibliography}
\end{document}